\documentclass[11pt]{article}

\title{Notes about extended real- and set-valued functions}

\author{Andreas H. Hamel\thanks{Yeshiva University New York, Department of Mathematical
Sciences, hamel@yu.edu}, Carola Schrage\thanks{Martin--Luther--University
Halle--Wittenberg, Institute for Mathematics, carola.schrage@mathematik.uni-halle.de}}

\date{}

\usepackage{amsmath, amssymb}

\newtheorem{theorem}{Theorem}

\newtheorem{remark}[theorem]{Remark}

\newtheorem{definition}[theorem]{Definition}
\newtheorem{proposition}[theorem]{Proposition}
\newtheorem{example}[theorem]{Example}

\numberwithin{equation}{section} \numberwithin{theorem}{section}

\newcommand{\of}[1]{\ensuremath{\left( #1 \right)}}

\newcommand{\abs}[1]{\ensuremath{\left| #1 \right|}}
\newcommand{\cb}[1]{\ensuremath{ \left\{ #1 \right\} }}
\newcommand{\sqb}[1]{\ensuremath{ \left[ #1 \right] }}

\newcommand{\bs}{\backslash}

\newcommand{\pend}{\hfill $\square$}
\newcommand{\eps}{\ensuremath{\varepsilon}}

\newcommand{\vp}{\ensuremath{\varphi}}

\newcommand{\R}{\mathrm{I\negthinspace R}}
\newcommand{\OLR}{\overline{\mathrm{I\negthinspace R}}}

\newcommand{\dom}{{\rm dom \,}}
\newcommand{\epi}{{\rm epi \,}}
\newcommand{\hypo}{{\rm hypo \,}}
\newcommand{\gr}{{\rm graph \,}}
\newcommand{\cl}{{\rm cl \,}}

\newcommand{\co}{{\rm co \,}}

\newcommand{\isum}{{+^{\negmedspace\centerdot\,}}}
\newcommand{\ssum}{{+_{\negmedspace\centerdot\,}}}
\newcommand{\idif}{{-^{\negmedspace\centerdot\,}}}
\newcommand{\sdif}{{-_{\negmedspace\centerdot\,}}}

\newcommand{\triup}{{\rm \vartriangle}}
\newcommand{\trido}{{\rm \triangledown}}

\newcommand{\cone}{{\rm cone\,}}

\begin{document}
\maketitle

\begin{abstract}
An order theoretic and algebraic framework for the extended real numbers is established
which includes extensions of the usual difference to expressions involving $-\infty$
and/or $+\infty$, so-called residuations. Based on this, definitions and results for
directional derivatives, subdifferentials and Legendre--Fenchel conjugates for extended
real-valued functions are given which admit to include the proper as well as the improper
case. For set-valued functions, scalar representation theorems and a new conjugation
theory are established. The common denominator is that the appropriate image spaces for
set-valued functions share fundamental structures with the extended real numbers: They
are order complete, residuated monoids with a multiplication by non-negative real
numbers.
\\[.2cm]
\\[.2cm]
\end{abstract}

\section{Motivation and bibliographical comments}

Without any doubts, the notion of an extended real-valued function turned out to be
extremely useful in variational analysis, optimization theory and beyond. On the one
hand, several operations like taking the directional derivative or the infimal
convolution, even performed on real-valued or proper functions, may lead to functions
which also attain the values $+\infty$ and/or $-\infty$, and it would be really awkward
to exclude such cases (see already \cite[p. 167]{IoffeTikhomirov79}). On the other hand,
the added element $+\infty$ admits the inclusion of constraints in a very elegant and
concise way (compare \cite[p. 23]{Rockafellar70}).

Almost all textbooks and relevant papers on convex and variational analysis make use of
this notion. As examples we mention \cite{Moreau66}, \cite{IoffeTikhomirov68},
\cite{Rockafellar70}, \cite{Rockafellar74}, \cite{EkelandTemam76},
\cite{IoffeTikhomirov79}, all already published before 1980.

To our opinion, the most thorough investigation of extended real-valued functions has
already been made by Jean Jacques Moreau in \cite{Moreau63a}, \cite{Moreau63b},
\cite{Moreau66}. It is a stunning and startling fact that his attempt to "algebraize" the
extended reals was not exploited consequently later on. Compare \cite[p.
6]{Rockafellar74}, and \cite[p. 9]{Aubin79}, where the  operation $\of{+\infty} +
\of{-\infty}$ is called "risky" and "undefined and forbidden", respectively, and also the
classic \cite[p. 24]{Rockafellar70}. Even in the more recent \cite[p.
15]{RockafellarWets98}, the authors state "there's no single, symmetric way of handling
$\infty - \infty$".

Most authors try to avoid the difficulties (like in \cite[p. 38]{Phelps93}: "we won't
have occasion to worry about $+\infty - \infty$ or $0\cdot\infty$") by restricting the
theory to proper functions or just ignore the problem. An extreme with respect to this
"avoiding approach" is the standard volume on infinite dimensional analysis \cite[p.
2]{AliprantisBorder06} which reads "The combination $+\infty - \infty$ of symbols has no
meaning. The symbols $+\infty$ and $-\infty$ are not really meant to be used for
arithmetic, they are only used to avoid awkward expressions involving infima and
suprema." In this note, we show that just the opposite works well.

To avoid the development (or the use) of an arithmetic for the extended reals does not
only passes a chance, it may also lead to imprecise statements. An example for the latter
can even be found in otherwise impressive textbooks: Theorem 2.3.1, (ix) in
\cite{Zalinescu02} does not hold for improper functions (no such assumption made in the
quoted reference) unless one uses the inf-addition on the left and the sup-addition on
the right hand side (see below for definitions) of the equation $\of{f \square g}^* = f^*
+ g^*$. The same remark applies, for example, to the first part of theorem 1 in
\cite[Section 3.4, p. 178]{IoffeTikhomirov79}.

In this note, we give an extension of Moreau's approach to extended real-valued functions
by noting that the correct algebraic framework is an order complete, residuated monoid
with a multiplication with non-negative real numbers. The advantage of this complicated
sounding construct is manifold: First, there is no need anymore to "explain away" the
value $-\infty$ or to introduce algebraic rules for expressions like $\of{+\infty} +
\of{-\infty}$ "by convention" (\cite[p. 15]{RockafellarWets98}), or to avoid them.
Secondly, new operations can be introduced which give a precise meaning to expressions
like $\of{+\infty} - \of{-\infty}$, and one obtains a whole calculus for addition and
residuation/difference in $\OLR = \R\cup\cb{-\infty}\cup\cb{+\infty}$. On an abstract
level, some of these observations have already been made by Mart{\'{\i}}nez-Legaz, Singer
and Getan in \cite{MartinezLegazSinger95} and \cite{GetanMartinezLegazSinger03}.

Moreover, our approach will also simplify the notation avoiding symbols like $r \dotplus
- s$ (see, among others, Moreau's papers, \cite[example 2.3]{MartinezLegazSinger95},
\cite{Singer97}). Finally, it will become clear that the theory is completely symmetric,
because our slightly different point of view (compared to Mart{\'{\i}}nez-Legaz, Singer
and others) is that there are two different ways for introducing algebraic and order
structures in $\OLR$ and more general sets as shown in the section about set-valued
functions. This follows Moreau's original idea of defining convex and concave functions
using different additions and image spaces.

Since we consequently work with two algebraically different copies of the extended reals,
we have to say which of the two is used as an image space if we define an extended
real-valued function. Thus, there are two classes of such functions. Not very
surprisingly, the multiplication by $-1$ transfers a function of one class into one of
the other, an operation which is nothing else than a duality in the sense of
\cite{Singer97}. We show how these concepts can be used, for example, to define
directional derivatives and subdifferentials of improper convex/concave functions in a
coherent way.

A second new feature of our approach is that in order to obtain complete (duality)
results the set of dual variables is extended by improper elements and, moreover the
definition of the Legendre-Fenchel conjugate is altered: the new definition involves an
additional real variable which comes from the idea that the conjugate should be defined
on the set of affine functions rather than on the set of linear functions. It does not
make a difference if the function is proper, but is does if not since the improper
"linear" functions are not additive. A somehow surprising result is that the conjugate of
the infimal convolution of two functions turns out to be the supremal convolution of
their conjugates -- with respect to the new primal variable.

We mention that improper affine functions have been used in \cite{Oettli82} in order to
formulate duality results for optimization problems involving set-valued maps. We are not
aware of further references, but we think there should be some.

Finally, we consider set-valued functions and give an extension of the theory formulated
in \cite{Hamel09} to improper set-valued functions using the (improper) scalar ones. In
fact, the present note has been written since we wanted to have a coherent framework for
proper and improper scalarizations of (closed convex) set-valued functions. The approach
follows ideas of \cite{Schrage09Diss}: In particular, using the representation of
set-valued closed convex functions by families of extended real-valued ones we give a new
definition of Legendre-Fenchel conjugates for set-valued functions and conclude with a
Fenchel-Moreau theorem which includes the proper as well as the improper case.

We conclude the introduction by noting that it does not take more than 5 pages and
relatively elementary mathematics to introduce the two possibilities for an algebraic and
order theoretic framework in $\OLR$, which seems affordable for classroom and textbook
purposes.

\section{A basic result from residuation theory}

In this section, we consider a lattice ordered set with an algebraic operation which we
call addition, denoted by $+$, since it corresponds to "usual" additions in most special
cases we have in mind. In the following, we understand by a partially ordered groupoid a
nonempty set $W$ with a binary relation $+ \colon W \times W \to W$ and a partial order
$\leq$ which are compatible: $u,v,w \in W$ and $u \leq v$ imply $u + w \leq v + w$. The
sum $u + M$ of $u \in W$ and $M \subseteq W$ is understood in the Minkowski sense with
$u+M = \emptyset$ if $M = \emptyset$. The following theorem can be extracted, for
example, from \cite[chapter XII]{Fuchs63}. See also \cite[proposition
2.6]{MartinezLegazSinger95}, \cite[proposition 2.1]{GetanMartinezLegazSinger03} for parts
(c), (d).

\begin{theorem}
\label{ThmInfRes} Let $\of{W, +, \leq}$ be a partially ordered commutative groupoid. The
following statements are equivalent:
\\
(a) For each $u, v \in W$ there is $w \in W$ such that for $w' \in W$
\[
u \leq v + w' \quad \Longleftrightarrow \quad w \leq w';
\]
(b) For each $u, v \in W$ the set $\cb{w' \in W \colon u \leq v + w'}$ has a least
element;
\\
(c) For $u \in W$ and $M \subseteq W$ such that $\inf M$ exists it holds
\[
u + \inf M = \inf\of{u + M};
\]
(d) For each $u, v \in W$ there exists $\inf\cb{w' \in W \colon u \leq v + w'} \in W$ and
it holds
\[
u \leq v + \inf\cb{w' \in W \colon u \leq v + w'}.
\]
\end{theorem}

{\sc Proof.} The equivalence of (a) and (b) is obvious. Assume (a). Then
\[
\forall m \in M \colon u + \inf M \leq u + m,
\]
hence $u + \inf M$ is a lower bound of $u + M$. On the other hand, let $w \in W$ such
that
\[
\forall m \in M \colon w \leq u + m.
\]
By assumption, there is $\bar{w} \in W$ such that
\[
w \leq u + w' \; \Leftrightarrow \; \bar{w} \leq w'.
\]
Hence $\bar{w} \leq m$ for all $m \in M$ and therefore $\bar{w} \leq \inf M$. Again by
assumption $w \leq u + \inf M$ which proves that $u + \inf M$ is the infimum of $u + M$.

Assuming (c) we define $w = \inf\cb{w' \in W \colon u \leq v + w'}$. Then
\[
u \leq \inf\cb{v + w' \colon w' \in W, \; u \leq v + w'} \underset{\text{(c)}}{=}
    v + \inf\cb{w' \in W \colon u \leq v + w'} = v + w.
\]

Now, from (c) it follows
\[
v \leq \inf\cb{u + w' \in W \colon v \leq u + w'}
    = u + \inf\cb{w' \in W \colon v \leq u + w'},
\]
i.e. (d) holds true.

Finally, assume (d) and define $M = \cb{w \in W \colon u \leq v + w}$. Then, by
assumption, $\inf M \in W$ and
\[
u \leq v + w \; \Leftrightarrow \; w \in M \; \Leftrightarrow \;
    \inf M \leq w
\]
which proves (a).

This completes the proof of the theorem. \pend

\medskip We shall call a partially ordered commutative groupoid satisfying the conditions
of theorem \ref{ThmInfRes} inf--residuated. The following theorem can be proven with
parallel arguments and gives conditions for sup--residuated groupoids.

\begin{theorem}
\label{ThmSupRes} Let $\of{W, +, \leq}$ be a partially ordered commutative groupoid. The
following statements are equivalent:
\\
(a) For each $u, v \in W$ there is $w \in W$ such that for $w' \in W$
\[
 v + w' \leq u \quad \Longleftrightarrow \quad w' \leq w;
\]
(b) For each $u, v \in W$ the set $\cb{w' \in W \colon v + w' \leq u}$ has a greatest
element;
\\
(c) For $u \in W$ and $M \subseteq W$ such that $\sup M$ exists it holds
\[
u + \sup M = \sup\of{u + M};
\]
(d) For each $u, v \in W$ there exists $\sup\cb{w' \in W \colon v + w' \leq u} \in W$ and
it holds
\[
v + \sup\cb{w' \in W \colon v + w' \leq u} \leq u.
\]
\end{theorem}

\section{An algebraic approach}

\subsection{Order extension}

Adding two elements $-\infty$, $+\infty$ to the set $\R$ of real numbers we consider the
set $\OLR = \R\cup\cb{-\infty}\cup\cb{+\infty}$ and extend the usual order relations
$\leq$, $<$ on $\R$ to $\OLR$ by setting

\begin{align*}
\forall r \in \OLR & \colon -\infty \leq r \leq +\infty \\
\forall r \in \R & \colon -\infty < r < +\infty.
\end{align*}
With this extension of $\leq$, $\of{\OLR, \leq}$ becomes a partially ordered, complete
lattice: Every subset has an infimum and a supremum. In particular,
\begin{align}
\inf \emptyset = \sup \R = \sup \OLR = +\infty, \\
\inf \OLR = \inf \R = \sup\emptyset = -\infty.
\end{align}
Note that the commonly used conventions $\inf \emptyset = +\infty$, $\sup\emptyset =
-\infty$ are the unavoidable choice if one wants to maintain the following monotonicity
property: $M \subseteq N \subseteq \OLR$ implies $\inf M \geq \inf N$ and $\sup M \leq
\sup N$.

\subsection{Addition}

There are two ways to extend the addition from $\R$ to $\OLR$ by means of the order
relation $\leq$. We obtain two different algebraic operations in $\OLR$.

\begin{definition}
\label{DefInfSubAdd} The binary operations $\isum \colon \OLR \times \OLR \to \OLR$ and
$\ssum \colon \OLR \times \OLR \to \OLR$ defined by
\begin{align}
r \isum s & = \inf\cb{a + b \colon a, b \in \R, \; r \leq a, \; s \leq b} \\
r \ssum s & = \sup\cb{a + b \colon a, b \in \R, \; a \leq r, \; b \leq s}
\end{align}
for $r, s \in \OLR$ are called the inf-addition and the sup-addition in $\OLR$,
respectively.
\end{definition}

The terminology is due to \cite{RockafellarWets98}. Already Moreau \cite{Moreau63a}
introduced the two different additions in $\OLR$. Clearly, both operations coincide with
the usual addition on $\R$. The notable differences are
\begin{align*}
\of{+\infty} \isum \of{-\infty} & = \of{-\infty} \isum \of{+\infty} = +\infty, \\
\of{+\infty} \ssum \of{-\infty} & = \of{-\infty} \ssum \of{+\infty} = -\infty.
\end{align*}
Since inf-adding $+\infty$ always gives $+\infty$ as a result, we say that $+\infty$
dominates the inf-addition. Likewise, $-\infty$ dominates the sup-addition. Both
operations are compatible with the order $\leq$ on $\OLR$ in the usual sense. Thus,
$\of{\OLR, \isum, \leq}$ and $\of{\OLR, \ssum, \leq}$ are ordered commutative monoids
which are complete lattices. The following result describes the relationships between
inf-/sup-addition and the order relation. Compare, for example, \cite{Moreau63a},
proposition 1 and 2 with $F=M$, $G = N$ and $f\of{x} = g\of{x}= x$.

\begin{proposition}
\label{ProSumOrder} Let $M, N \subseteq \OLR$. Then
\begin{align}
\label{EqInfSupIPlus} \inf \of{M \isum N} = \inf M \isum \inf N, \quad
    & \sup \of{M \isum N} \leq \sup M \isum \sup N  \\
\label{EqInfSupSPlus} \inf M \ssum \inf N \leq \inf \of{M \ssum N}, \quad
    & \sup M \ssum \sup N = \sup \of{M \ssum N}
\end{align}
where the sum of sets is understood in the Minkowski sense.
\end{proposition}

{\sc Proof.} By definition,
\[
\inf M \isum \inf N=\inf\cb{r+s:\; r,s\in\R,\, \inf M\leq r, \inf N\leq s}
\]
and for all $m\in M, n\in N$ holds $\inf M \isum \inf N\leq m\isum n$, thus $\inf M \isum
\inf N \leq \inf\of{M \isum N}$. On the other hand, if $M = \emptyset$ or $N =
\emptyset$, then both of $\inf M \isum \inf N$ and $\inf\of{M \isum N}$ give the result
$+\infty$ since the latter element dominates the inf-addition. If $M, N \neq \emptyset$
and one of them has infimum $-\infty$ then $\inf M \isum \inf N = \inf\of{M \isum N} =
-\infty$. If $-\infty < \inf M, \inf N$ then, for each $\eps > 0$ there are $m_\eps \in
M\cap \R$, $n_\eps \in N \cap \R$ such that $m_\eps + n_\eps \leq \inf M \isum \inf N +
\eps$, hence $\inf M \isum \inf N \leq \inf\of{M \isum N}$ which gives, together with the
first part, equality. \pend

\medskip In particular, with $r \in \OLR$, $N = \cb{r}$ we obtain (again, compare
\cite{Moreau63a}, p. 7, formulas (2.9), (2.12))
\begin{equation}
\label{EqConstPlusInf} r \isum \inf M = \inf\of{\cb{r} \isum M}, \quad r \ssum \sup M =
\sup\of{\cb{r} \ssum M}.
\end{equation}
Finally, note that the inequalities in \eqref{EqInfSupIPlus} and \eqref{EqInfSupSPlus}
are not satisfied as equations in general. A counterexamples can already be found in
\cite{Moreau63b}, p. 7.

\subsection{Multiplication with $-1$}

By setting
\[
\of{-1}\of{+\infty} = -\infty, \quad \of{-1}\of{-\infty} = +\infty
\]
we extend the multiplication of real numbers with $-1$ to $\OLR$. As usual, we abbreviate
$\of{-1}r$ to $-r$ for $r \in \OLR$ if no confusion arises. Obviously,
\[
\forall r, s \in \OLR \colon \of{r \leq s \; \Leftrightarrow \; -s \leq -r}
\]
and hence for each $M \subseteq \OLR$
\begin{equation}
\label{EqMinusInfSup} \of{-1} \inf M = \sup \of{-1}M.
\end{equation}
Thus, the multiplication with $-1$ is a duality of $\OLR$ onto itself in the sense of
\cite{Singer97}, chapter 5. Here and in the following, we make use of $\of{-1}M = \cb{-r
\colon r \in M}$ if $M \neq \emptyset$ and $\of{-1}M = \emptyset$ if $M = \emptyset$.

\begin{proposition}
\label{PropMinusSum} For $r,s \in \OLR$,
\[
\of{-1}\of{r \isum s} = \of{-1}r \ssum \of{-1}s.
\]
\end{proposition}

{\sc Proof.} Using definition \ref{DefInfSubAdd} and \ref{EqMinusInfSup} above we obtain
\begin{align*}
\of{-1}\of{r \isum s} & = \of{-1}\inf\cb{a + b \colon a, b \in \R, \; r \leq a, \; s \leq b} \\
    & = \sup\cb{\of{-a} + \of{-b} \colon a, b \in \R, \; r \leq a, \; s \leq b} \\
    & = \sup\cb{a' + b' \colon a', b' \in \R, \; a' \leq -r, \; b' \leq -s} \\
    & = \of{-1}r \ssum \of{-1}s
\end{align*}
which already proves the claim. \pend

\subsection{Residuation and inf-/sup-difference}

Proposition \ref{ProSumOrder} together with theorem \ref{ThmInfRes} and \ref{ThmSupRes}
tell us that $\of{\OLR, \isum, \leq}$ and $\of{\OLR, \ssum, \leq}$ are residuated
semigroups with a neutral element (i.e. residuated monoids), see e.g. \cite{Fuchs63},
chap. XII. The corresponding residuation operations may serve as extensions of the
difference from $\R$ to $\of{\OLR, \isum, \leq}$ and $\of{\OLR, \ssum, \leq}$,
respectively. This motivates the following definition. Residuation operations for $x
\mapsto r + x$ in $\OLR$ have not been considered by Moreau. Only Mart{\'{\i}}nez-Legaz,
Singer and Getan (see \cite{MartinezLegazSinger95}, \cite{GetanMartinezLegazSinger03})
seem to have realized the importance of residuation for the foundation of convex
analysis.

\begin{definition}
\label{DefInfSupDifference} Let $r, s \in \OLR$. The element
\[
r \idif s = \min\cb{t \in \OLR \colon r \leq s \isum t}
\]
is called the inf-difference of $r$ and $s$. The element
\[
r \sdif s = \max\cb{t \in \OLR \colon s \ssum t \leq r}
\]
is called the sup-difference of $r$ and $s$.
\end{definition}

One easily obtains for all $r,s \in \OLR$
\begin{align}\label{EQAlternativeDifference}
r \idif s &= \inf\cb{t \in \R \colon r \leq s \isum t}, \\
r \sdif s &= \sup\cb{t \in \R \colon s \ssum t \leq r}
\end{align}
and
\begin{align}
\label{EqMinusInftyI}
 r \idif \of{+\infty} = -\infty, \quad & \of{-\infty} \idif r = -\infty, \\
\label{EqMinusInftyS}
 r \sdif \of{-\infty} = +\infty, \quad & \of{+\infty} \sdif r = +\infty.
\end{align}
The rules for a subtraction of least and greatest elements from each other are as shown
below:
\begin{align*}
   \of{+\infty} \idif \of{-\infty} = +\infty, & \quad \of{+\infty} \sdif \of{+\infty} = +\infty, \\
   \of{+\infty} \idif \of{+\infty} = -\infty, & \quad \of{+\infty} \sdif \of{-\infty} = +\infty, \\
   \of{-\infty} \idif \of{+\infty} = -\infty, & \quad \of{-\infty} \sdif \of{-\infty} = +\infty, \\
   \of{-\infty} \idif \of{-\infty} = -\infty, & \quad \of{-\infty} \sdif \of{+\infty} =
   -\infty.
\end{align*}

Moreover, from theorem \ref{ThmInfRes} (a) with $u = r \in \OLR$, $v = s \in \OLR$, $w =
r \idif s$ and $w' = 0$, and likewise with the help of theorem \ref{ThmSupRes}, we obtain
\begin{equation}
\label{EqIdifInequ}
 r \leq s \; \Leftrightarrow \; r \idif s \leq 0 \; \Leftrightarrow \; 0 \leq s \sdif r.
\end{equation}

The following result gives relationships between  inf-/sup-addition, inf-/sup-subtraction
and multiplication with $-1$.

\begin{proposition}
\label{ProDualDiff} Let $r, s \in \OLR$. Then
\begin{align}
\label{EqProDualDiff1}
r \idif s & = r \ssum \of{-1}s, \\
\label{EqProDualDiff2} r \sdif s & = r \isum \of{-1}s, \\
\label{EqProDualDiff3} s \sdif r & = (-1)\of{r \idif s}, \\
s \idif r & = (-1)r \idif (-1)s, \\
s \sdif r & = (-1)r \sdif (-1)s.
\end{align}
\end{proposition}

{\sc Proof.} If $r=-\infty$, or if $s=+\infty$, then $r \idif s  = r \ssum
\of{-1}s=-\infty$, see \eqref{EqMinusInftyI}. If $r=+\infty$ and $s<+\infty$, or if
$s=-\infty$ and $r>-\infty$, then $r \idif s  = r \ssum \of{-1}s=+\infty$. This proves
\eqref{EqProDualDiff1} since if $r,s\in \R$, the formula is known to be true. Likewise,
\eqref{EqProDualDiff2} is proven.

Next, we use \eqref{EqMinusInfSup} and \eqref{EQAlternativeDifference} to obtain
\[
(-1)(r\idif s)=\sup\cb{t\in \R \colon r \leq s\isum (-1)t}.
\]
Since it suffices to take the supremum over $t \in \R$ we get
\[
r \leq s\isum (-1)t \quad \Leftrightarrow \quad r \ssum t \leq s
\]
and \eqref{EqProDualDiff3} follows from definition \ref{DefInfSupDifference}.

The last two equations are immediate from \eqref{EqProDualDiff1}, \eqref{EqProDualDiff2}.
\pend

\medskip We establish a calculus for manipulating inf-/sup-differences. Since these
operations are special cases of residuation mappings these rules are well-known, see for
example \cite[Lemma 3.2]{BlountTsinakis03} where the sup versions can be found.

\begin{proposition}
\label{ProDiffCalc1} (a) For each $r\in\OLR$,
\[
r\idif r = \left\{
    \begin{array}{ccc}
    0 & : & r \in \R \\
    -\infty & : & r \not\in \R
    \end{array}
    \right., \quad
r\sdif r = \left\{
    \begin{array}{ccc}
    0 & : & r \in \R \\
    +\infty & : & r \not\in \R
    \end{array}
    \right. .
\]
(b) For each $r, s, t \in \OLR$ with $r\leq s$,
\[
r \idif t \leq s \idif t, \quad t \idif s \leq t \idif r,
\]
and
\[
r \sdif t \leq s \sdif t, \quad t\sdif s \leq t\sdif r.
\]
(c) For each $a,b,r,s \in\OLR$,
\[
(a\isum r) \idif (b\isum s) \leq (a\idif b) \isum (r \idif s),
\]
especially
\[
(a\isum r) \idif (r\isum b) \leq (a\idif b).
\]
Symmetrically,
\[
(a\sdif b) \ssum (r \sdif s) \leq (a\ssum r) \sdif (b\ssum s),
\]
especially
\[
(a\sdif b) \leq (a\ssum r) \sdif (r\ssum b).
\]
(d) If $M\subseteq \OLR$ and $a\in \OLR$, then
\begin{align*}
a\idif \inf M & = \sup\limits_{m\in M} (a\idif m), \\
a\sdif \sup M & = \inf\limits_{m\in M} (a\sdif m).
\end{align*}
\end{proposition}

{\sc Proof.} (a) is straightforward from the definitions. For (b) observe that if $r \leq
s$, then $s \leq t \isum t'$ implies $r \leq t \isum t'$, and this in turn gives $r \idif
t \leq s \idif t$. Moreover, if $r \leq s$, then $t \leq r \isum t'$ implies $t \leq s
\isum t'$, and this in turn gives $t \idif s \leq t \idif r$. The relationships for
$\sdif$ can be proven similarly. We turn to (c). Take $t_1, t_2 \in \OLR$ such that $a
\leq b \isum t_1$ and $r \leq s \isum t_2$. Then $\of{a \isum r} \leq \of{b \isum s}
\isum \of{t_1+t_2}$, hence
\[
\of{a \isum r} \idif \of{b \isum s} \leq t_1 \isum t_2.
\]
Taking the infimum over $t_1$ satisfying $a \leq b \isum t_1$ and $t_2$ satisfying $r
\leq s \isum t_2$ gives the result. The second formula in (c) is immediate by setting
$r=s$ in the first and applying (a). The results for $\sdif$ are proven likewise. (d)
Using \eqref{EqProDualDiff1}, \eqref{EqProDualDiff2} as well as \eqref{EqMinusInfSup} and
(e) of proposition \ref{ProDualDiff} we obtain
\[
a \idif \inf M = a \ssum \of{-1}\inf M = a \ssum \sup -M =
    \sup\of{a \ssum \of{-M}} = \sup\of{a \idif M}.
\]
The second equation follows from the first since by \eqref{EqProDualDiff2} we have $a
\sdif \sup M = a \isum \inf (-1)M$ and thus $a \sdif \sup M= \inf\limits_{m\in M}\of{a
\sdif m}$. \pend

\subsection{Multiplication with non-negative reals}

The multiplication with non-negative reals is extended to $\OLR$ by
\[
\forall t>0 \colon t \cdot \of{\pm\infty} = \pm\infty
\]
and $0 \cdot \of{\pm\infty} = 0 \in \R$. The triples $\of{\OLR, \isum, \cdot}$ and
$\of{\OLR, \ssum, \cdot}$ are conlinear spaces (see appendix for a definition) consisting
only of convex elements. That is, the multiplication with non-negative real numbers
distributes over $\isum$ and $\ssum$ as well as the other way round.

The order relation $\leq$ as defined above is compatible with this algebraic structure in
the usual sense. We write $\R^\vartriangle = \of{\OLR, \isum, \cdot, \leq}$,
$\R^\triangledown = \of{\OLR, \ssum, \cdot, \leq}$, and drop the $\cdot$ for
multiplication if no confusion arises.

\begin{proposition}
\label{ProDiffCalc2} For all $a, b\in\OLR$ and $t\geq 0$ it holds
\[
t(a\idif b)=ta\idif tb, \quad t(a\sdif b)=ta\sdif tb.
\]
\end{proposition}

{\sc Proof.} The relationships are trivial for $t = 0$. If $t > 0$ then
\begin{align*}
t(a\idif b) & = \inf \cb{ts \in \R \colon a \leq b\isum s} \\
            & = \inf \cb{s' \in \R \colon ta \leq tb \isum s'} = ta \idif tb.
\end{align*}
The result for $\sdif$ follows similarly. \pend

\section{Extended real-valued functions}

From the above, it should be clear that there are two types of extended real-valued
functions, those mapping into $\R^\vartriangle$ and those mapping into
$\R^\triangledown$. The point-wise multiplication with $-1$ transfers a function of one
class into a function of the other. This point of view differs slightly from
\cite{GetanMartinezLegazSinger03} where (only one copy of) $\OLR$ is considered with two
additions and two corresponding residuation operations. This might appear to be just a
tiny shift of weight, but it becomes important when it comes to set-valued functions: The
replacements of $\R^\triup$ and $\R^\trido$ will have rather different looking elements.

\begin{definition}
\label{DefDomEpi} Let $f \colon X \to \OLR$. The epigraph and the hypographs of $f$ are
the sets
\[
\epi f = \cb{\of{x, r} \in X \times \R \colon f\of{x} \leq r}, \quad
    \hypo f = \cb{\of{x, r} \in X \times \R \colon r \leq f\of{x}},
\]
respectively. The effective domain of a function $g \colon X \to \R^\vartriangle$ is the
set
\[
\dom g = \cb{x \in X \colon g\of{x} < +\infty}
\]
whereas the effective domain of a function $h \colon X \to \R^\triangledown$ is the set
\[
\dom h = \cb{x \in X \colon -\infty < h\of{x}}.
\]
\end{definition}

The concept of the domain depends on the image space, so if one knows the latter, one
also knows which definition to use. Therefore, we do not introduce different symbols.
Note also that the collection of functions into $\R^\vartriangle$ ($\R^\triangledown$) is
a conlinear space under point-wise addition $\isum$ ($\ssum$) and multiplication with
non-negative reals, but neither collection is a linear space. Mixing up the image spaces
may lead to strange effects as the next example shows.

\begin{example}
Consider the functions $f, g \colon \R \to \R^\vartriangle$ defined by
\[
f\of{x} = \left\{
    \begin{array}{lll}
    +\infty & : & x < 2 \\
    -\infty & : & x \geq 2
    \end{array}
    \right., \quad
g\of{x} = \left\{
    \begin{array}{lll}
    +\infty & : & \abs{x} \leq 1 \\
    -\infty & : & \abs{x} > 1
    \end{array}
    \right.
\]
Both functions have a convex epigraph. So has the function $x \mapsto f\of{x} \isum
g\of{x}$. However, the function $x \mapsto f\of{x} \ssum g\of{x}$ neither has a convex
epigraph nor a convex hypograph.
\end{example}

The previous example also shows that convexity and $\isum$ are linked as are concavity
and $\ssum$. This justifies the following definition (already in \cite[p. 19]{Moreau63a}.

\begin{definition}
\label{DefConvexFunction} A function $g \colon X \to \R^\vartriangle$ is called convex if
\[
\forall t \in \of{0,1}, \; \forall x_1, x_2 \in X \colon
    g\of{tx_1 + \of{1-t}x_2} \leq tg\of{x_1} \isum \of{1-t}g\of{x_2}.
\]
A function $h \colon X \to \R^\triangledown$ is called concave if
\[
\forall t \in \of{0,1}, \; \forall x_1, x_2 \in X \colon
    th\of{x_1} \ssum \of{1-t}h\of{x_2} \leq h\of{tx_1 + \of{1-t}x_2}.
\]
\end{definition}

Again, the collection of all convex functions into $\R^\vartriangle$ (concave functions
into $\R^\triangledown$) is a conlinear space under point-wise addition $\isum$ ($\ssum$)
and multiplication with non-negative reals. Apparently, it does not make much sense to
consider convex functions into $\R^\triangledown$ or concave into $\R^\vartriangle$.

A function $g \colon X \to \R^\triup$ is convex if and only if $-g \colon X \to
\R^\trido$ is concave. A function $g$ mapping into $\R^\triup$ is convex if and only if
its epigraph is a convex subset of the linear space $X \times \R$, and a function $h$
mapping into $\R^\trido$ is concave if and only if its hypograph is convex.

\begin{definition}
\label{DefSublinFunction} A function $f \colon X \to \OLR$ is called positively
homogeneous if
\[
\forall t >0, \; \forall x \in X \colon
    f\of{tx} = tf\of{x}.
\]
A positively homogeneous convex function into $\R^\triup$ is called sublinear, and a
positively homogeneous concave function into $\R^\trido$  is called superlinear.
\end{definition}

Notice that we do not include the case $t=0$ in the definition of positive homogeneity.
Thus,
\[
g(x)=\left\{
       \begin{array}{ccc}
         -\infty & : & x \leq 0, \\
         +\infty & : & x > 0
       \end{array}
     \right.
\]
is a positively homogeneous function, while $0\cdot g(x)\neq g(0\cdot x)$ holds for all
$x\in X$.

\begin{example} {\bf (improper affine functions)}
\label{ExImpropAffine} Let $X$ be a topological linear space and $X^*$ its topological
dual. We write $x^*\of{x}$ for the value of an element $x^* \in X^*$ at $x \in X$. Let $r
\in \R$ and set $x^*_r(x)=x^*(x)-r$ for $x\in X$. Each $x^* \in X^*$ generates a closed
improper function $\hat x^*_r \colon X \to \R^\triup$ given by
\[
\hat x_r ^*\of{x} = \left\{
       \begin{array}{lll}
         -\infty & : & x_r^*\of{x} \leq 0 \\
         +\infty & : & x_r^*\of{x} > 0
       \end{array}
     \right.
\]
which we call inf-extension of the affine function $x \mapsto x_r^*\of{x}$. Analogously,
the improper sup-extension of $x \mapsto x_r^*\of{x}$ (with a closed hypograph) mapping
into $\R^\trido$ can be obtained by reversing the roles of $-\infty$ and $+\infty$. The
functions $\hat x^*_0$ are positively homogeneous, subadditive and superadditive, but not
additive, i.e. in general $\hat x^*_0\of{x_1 + x_2} \neq \hat x^*_0\of{x_1} \isum \hat
x^*_0\of{x_2}$ for $x_1, x_2 \in X$. Below, this will force us to define
Legendre--Fenchel conjugates acting on the set of affine rather than linear functions.
\end{example}

In the sequel, we shall write $x^*$ for $x^*_0$ and $\hat x^*$ for $\hat x^*_0$. Define
\[
\hat X^* = \cb{\hat x^* \colon x^* \in X^*}
\]
and $X^\triup = X^* \cup \hat X^*$. The set $X^\triup$ is called the (topological)
inf-dual of $X$ (correspondingly, the sup-dual $X^\trido$ of $X$ can be defined using the
sup-extensions of continuous linear functions). On $X^\triup$, an addition can be
introduced by
\[
\xi + \eta = \eta + \xi =
    \left\{
    \begin{array}{ccc}
    \hbox{inf-extension of} \; x^*+y^* & : & \xi=\hat x^*, \eta = \hat y^*, \\
    \xi & : & \xi=\hat x^*, \eta = y^*, \\
    x^* + y^* & : & \xi = x^*, \eta = y^*
    \end{array}
    \right.
\]
for $\xi, \eta \in X^\triup$ with $x^*, y^* \in X^*$. Taking the multiplication by
non-negative numbers pointwise, $\of{X^\triup, +, \cdot}$ is a conlinear space with
neutral element $0 \in X^*$.

The following representation formulas for (proper and improper) affine functions will be
used later on.

\begin{proposition} Let $\xi \in X^\triup$, $r \in \R$. Then
\begin{align}
\label{EqAffinSupConv1} \forall x_1,x_2\in X \colon
    \xi_r\of{x_1+x_2} & = \sup\limits_{r_1+r_2=r} \sqb{\xi_{r_1}\of{x_1} \ssum \xi_{r_2}\of{x_2}}, \\
\label{EqAffinSupConv2} \forall x_1,x_2\in X \colon
    \xi_r\of{x_1-x_2} & = \sup\limits_{r_1+r_2=r} \sqb{\xi_{r_1}\of{x_1} \idif
    \xi_{-r_2}\of{x_2}}.
\end{align}
\end{proposition}

{\sc Proof.} Exemplarily, we prove \eqref{EqAffinSupConv1}. The formula is obvious for
$\xi = x^* \in X^*$ and for $\xi = \hat 0^*$. Let $\xi = \hat x^* \in \hat X^*\bs\cb{\hat
0^*}$. Then, by definition of $\hat x^*$ and since $-\infty$ dominates the sup-addition
$\sup\cb{\hat x_{r_1}^*\of{x} \ssum \hat x_{r_2}^*\of{y} \colon r_1+r_2= r} = +\infty$ if
and only if there is $r_1+r_2\in \R$ such that $x^*\of{x}-r > 0$ and $x^*\of{y}-r > 0$
and $r_1+r_2=r$. But this is equivalent to $x^*\of{x+y}-r > 0$, and this in turn to $\hat
x_r^*\of{x+y} = +\infty$. Otherwise, $\hat x_r^*\of{x+y} =
\sup\limits_{r_1+r_2=r}\of{\hat x_{r_1}^*\of{x} \ssum \hat x_{r_2}^*\of{y}} = -\infty$.
\pend

\medskip Replacing sup by inf, $\ssum$ by $\isum$ and $\idif$ by $\sdif$ one may obtain
another pair of representation formulas.

\section{Applications}

\subsection{Directional derivatives of convex functions}

In this section, we shall show that the residuations $\idif$ and $\sdif$ may serve as
substitutes for the usual difference in the definition of the directional derivative of
extended real-valued convex and concave functions. Besides the obvious result (may be
useful or not) of having a coherent definition of directional derivatives also for
improper functions the constructions below indicate the deep connection between the order
relation in the image space, crucial for the definition of residuations, and directional
derivatives.

\begin{definition}
\label{DefDirDer} The directional derivative of a function $g \colon X \to
\R^\vartriangle$ at $x_0 \in X$ in direction $x \in X$ is given by
\[
g'\of{x_0, x} =
    \lim_{t \downarrow 0} \frac{1}{t}\sqb{\of{g\of{x_0 + tx} \idif g\of{x_0}}}.
\]
The directional derivative of a function $h \colon X \to \R^\triangledown$ at $x_0 \in X$
in direction $x \in X$ is given by
\[
h'\of{x_0, x} =
    \lim_{t \downarrow 0} \frac{1}{t}\sqb{\of{h\of{x_0 + tx} \sdif h\of{x_0}}}.
\]
\end{definition}

If $g\of{x_0} = +\infty$ then $g'\of{x_0, x} = -\infty$ since $r \idif \of{+\infty} =
-\infty$ for all $r \in \OLR$. If $g$ is convex and $g\of{x_0} = -\infty$ then
$g'\of{x_0, x} = -\infty$ iff there is $t>0$ such that $g\of{x_0 + tx} = -\infty$ and
$g'\of{x_0, x} = +\infty$ otherwise. If $g\of{x_0} \in \R$ the directional derivative
coincides with the classical one (see e.g. \cite{IoffeTikhomirov79}, p. 193). Similar
remarks apply to $h'\of{x_0, x}$.

\begin{remark}
We have $g'\of{x_0, x} = \of{-1}\of{-g}'\of{x_0, x}$ and similar for $h$. This can be
seen with the help of proposition \ref{ProDualDiff}.
\end{remark}

\begin{proposition}
\label{ProDirDer} If $g \colon X \to \R^\vartriangle$ is convex, then $x \mapsto
g'\of{x_0, x}$ is a sublinear function from $X$ into $\R^\vartriangle$. In this case,
\begin{equation}
\label{EqDirDerInf} g'\of{x_0, x} = \inf_{t>0}\frac{1}{t}\sqb{\of{g\of{x_0 + tx} \idif
g\of{x_0}}}.
\end{equation}
Likewise, if $h \colon X \to \R^\triangledown$ is concave, then $x \mapsto h'\of{x_0, x}$
is a superlinear function from $X$ into $\R^\triangledown$ and
\begin{equation}
h'\of{x_0, x} = \sup_{t>0}\frac{1}{t}\sqb{\of{h\of{x_0 + tx} \sdif h\of{x_0}}}.
\end{equation}
\end{proposition}

{\sc Proof.} The basic fact is the monotonicity of the difference quotient
\[
t \mapsto \frac{1}{t}\sqb{\of{g\of{x_0 + tx} \idif g\of{x_0}}}.
\]
Indeed, taking $t>0$ and $\tau \in \of{0,1}$ we obtain with the help of proposition
\ref{ProDiffCalc1}, (c) and proposition \ref{ProDiffCalc2}
\begin{align*}
\frac{1}{\tau t}\sqb{g\of{x_0 + \tau tx} \idif g\of{x_0}}
& \leq \frac{\tau}{\tau t}\sqb{g\of{x_0 + tx} \idif g\of{x_0}}\isum \frac{1-\tau}{\tau t}\sqb{\of{g(x_0)\idif g(x_0)}}\\
& \leq \frac{1}{t}\sqb{g\of{x_0 + tx} \idif g\of{x_0}}.
\end{align*}
Therefore, for $0 < s = \tau t < t$
\[
\frac{1}{s}\sqb{\of{g\of{x_0 + sx} \idif g\of{x_0}}} \leq
    \frac{1}{t}\sqb{\of{g\of{x_0 + tx} \idif g\of{x_0}}}
\]
and hence
\[
g'\of{x_0, x} = \inf_{t>0}\frac{1}{t}\sqb{\of{g\of{x_0 + tx} \idif g\of{x_0}}}
\]
holds true. For $s>0$,
\[
g'\of{x_0, sx} = s\cdot\inf_{st>0}\frac{1}{st}\sqb{\of{g\of{x_0 + stx} \idif g\of{x_0}}}
    = s g'\of{x_0, x},
\]
thus the directional derivative is positively homogeneous.

Finally, consider the function $g_t \colon X \to \R^\triup$ defined by
\[
g_t(x)= \frac{1}{t}\of{g(x_0 + tx) \idif g(x_0)} \; \text{for} \; x \in X.
\]
The epigraph of each function $g_t$ is convex and by the monotonicity of the difference
quotient,
\[
\epi g'(x_0,\cdot)=\bigcup\limits_{t>0}\epi g_t
\]
is convex. Thus, the directional derivative is a sublinear function.

The concave case is immediate, considering that $h$ is concave iff $-h$ is convex. \pend

\medskip The first part of the previous proposition is an extension of results like theorem 2.1.13 in
\cite{Zalinescu02}. In the "proper" theory, it is well-established that the linear
minorants of $x \mapsto g'\of{x_0,x}$ are precisely the elements of the subdifferential
of $g$ at $x_0$. Admitting improper affine functions (see example \ref{ExImpropAffine})
it is now possible to formulate an "improper" supplement for the "proper" theory. Recall
$X^\triup = X^* \cup \hat{X}^*$ (see example \ref{ExImpropAffine}).

\begin{definition} An element $\xi \in X^\triup$ is called an extended subgradient of
$g \colon X \to \R^\triup$ at $x_0 \in X$ iff
\begin{align}
\label{EqImpSubgrad1}
 \forall x \in X \colon \xi(x-x_0) \leq g\of{x} \idif g\of{x_0} 
\end{align}
The set of all subgradients of $g$ at $x_0$ is denoted by $\partial^{ex}g\of{x_0}$, the
set of improper subgradients of $g$ at $x_0$ is denoted by $\partial^{ip}g\of{x_0}$, and
the (classical) subdifferential is $\partial g\of{x_0} = \partial^{ex}g\of{x_0} \bs
\partial^{ip}g\of{x_0}$.
\end{definition}

Obviously, $-\infty \in \partial^{ex}g\of{x_0}$ for all $g$ and all $x_0 \in X$.


\begin{proposition}
\label{PropImpSubDer} Let $g \colon X \to \R^\triup$ be a convex function. The following
two statements are equivalent for $\xi \in X^\triup$, $x_0 \in X$:
\\
(a) $\xi$ is a subgradient of $g$ at $x_0$;
\\
(b) $\forall x \in X$: $\xi\of{x} \leq g'\of{x_0, x}$.
\end{proposition}

{\sc Proof.} First, assume (a) and choose $x = x_0 + ty$ with $t>0$, $y \in X$. We get
\[
\forall y \in X, \; \forall t>0 \colon \xi\of{ty} \leq g\of{x_0 + ty} \idif g\of{x_0}.
\]
Since $\frac{1}{t}\xi\of{ty} = \xi\of{y}$ for all $t>0$, $y \in X$ (no matter if $\xi$ is
proper or improper) we may conclude
\[
\forall y \in X, \; \forall t>0 \colon \xi\of{y} \leq \frac{1}{t}\sqb{g\of{x_0 + ty}
\idif g\of{x_0}}.
\]
Formula \eqref{EqDirDerInf} produces (b). Next, assume (b). Using \eqref{EqDirDerInf} and
choosing $t=1$ gives \eqref{EqImpSubgrad1}. \pend

\subsection{Legendre--Fenchel conjugation}

The following concepts and results supplement the conjugation theory for functions with a
proper closure. Throughout this section, we assume that $X$ is a separated locally convex
space with topological dual $X^*$. It is well-known that if $g$ is convex and improper,
then $\of{\cl g}\of{x} \in \cb{\pm\infty}$ for all $x \in X$ (see \cite[proposition
2.2.5]{Zalinescu02}). In this case, $\dom\of{\cl g} = \cl\of{\dom g}$.

\begin{theorem}
\label{ThmImpropAffMin} Let $g \colon X \to \R^\triup$ be improper closed convex. Then
$g$ is the pointwise supremum of its improper closed minorants, and if $g \not\equiv
-\infty$, then there are $x^* \in X^*\bs\cb{0}$ and $r \in \R$ such that
\[
\forall x \in X \colon \hat x_r^*\of{x} \leq g\of{x}.
\]
\end{theorem}

{\sc Proof.} If $g$ is $g \equiv +\infty$ or $\equiv -\infty$, then $g(x)=\hat0_{-1}^*$
or $g(x)=\hat 0^*$ for all $x\in X$, respectively. Thus, let us assume $\dom g \not\equiv
\emptyset$ and $g \not\equiv -\infty$. Since $\dom g$ is a closed convex subset of $X$ it
is the intersection of all closed half spaces including it. Each such half space has the
form
\begin{equation}
\label{EqDomHalfSpace} \cb{x \in X \colon x^*\of{x}-r \leq 0} \supseteq \dom g
\end{equation}
for some $x^* \in X^*$, $r \in \R$. The function $\hat x_r^*:X\to\OLR$ certainly is a
closed improper affine minorant of $g$. Since $\dom g$ is the intersection of all $\dom
x_r^*$ with $x^*$, $r$ satisfying \eqref{EqDomHalfSpace} the result follows. \pend

\medskip The next result characterizes improper affine minorants of improper
$\R^\triup$-valued functions.

\begin{theorem}
\label{ThmImpAffMin} Let $g \colon X \to \R^\triup$. The following statements are
equivalent for $\xi \in X^\triup$ and $r \in \R$:
\\
(a) $\forall x \in X$: $\xi_r\of{x} \leq g\of{x}$,
\\
(b) $\sup_{x \in X} \cb{\xi_r\of{x} \idif g\of{x}} \leq 0$,
\\
(c) $\sup_{x \in X} \cb{\hat \xi_r\of{x} \ssum \of{-1}g\of{x}} \leq 0$,
\\
(d) $\inf_{x \in X} \cb{g\of{x} \idif \xi_r\of{x}} \geq 0$,
\\
(e) $\inf_{x \in X} \cb{g\of{x} \ssum \of{-1}\xi_r\of{x}} \geq 0$.

Moreover, if $\xi = \hat x^*$ with $x^* \in X^*$ then (a) through (e) are equivalent to
\\
(f) $\dom \hat x_r^*\supseteq \dom g$,
\\
and the suprema in (b), (c) are $-\infty$ whereas the infima in (d), (e) are $+\infty$.
\end{theorem}

{\sc Proof.} The equivalence of (a) and (b) follows from \eqref{EqIdifInequ}.  The
equivalence of (b) and (c) is immediate from \eqref{EqProDualDiff1}. Likewise, (d) and
(e) are equivalent. The equivalence of (b) and (d) follows from \eqref{EqProDualDiff1}
and \eqref{EqMinusInfSup}.

If $\xi = \hat x^*$ then the difference $\hat x_r^*\of{x} \idif g\of{x}$ is $+\infty$ if
and only if $x \in \dom g\bs\dom\hat x_r^*$, so (b) and (f) are equivalent. \pend

\medskip The previous theorem together with the effect described in example \ref{ExImpropAffine}
gives rise to the following definition.

\begin{definition}
\label{DefConjugate} The Legendre-Fenchel conjugate of $g \colon X \to \R^\triup$ is the
function $g^* \colon X^\triup \times \R \to \R^\trido$ which is given by
\[
g^*\of{\xi, r} = \sup_{x \in X} \cb{\xi_r\of{x} \idif g\of{x}}.
\]
for $\xi \in X^\triup$, $r \in \R$.
\end{definition}

Take $x^* \in X^*$ and $r \in \R$. Of course,
\[
g^*\of{x^*, r} = \sup_{x \in X} \cb{x_r^*\of{x} \idif g\of{x}}
    = \sup_{x \in X} \cb{x^*\of{x} \idif g\of{x}} - r = g^*\of{x^*} - r
\]
where $g^*\of{x^*}$ is the classical Legendre-Fenchel conjugate of $g$ at $x^* \in X^*$.
Moreover (compare the previous theorem),
\begin{align}\label{EqImpropConjug}
g^*\of{\hat x^*, r} = \left\{
                   \begin{array}{ccc}
                   -\infty & : & \dom \hat x_r ^*\of{\cdot } \supseteq \dom g, \\
                   +\infty & : & \hbox{otherwise.}
                   \end{array}
                   \right.
\end{align}
Therefore, $\xi_r$ with $\xi \in X^\triup$, $r\in \R$ is a (proper or improper) closed
affine minorant of $g$ if and only if $g^*\of{\xi, r} \leq 0$.

\begin{proposition}\label{YoungFenchel}
Let $g \colon X \to \R^\triup$ be a function, $\of{\xi, r} \in X^\triup \times \R$. The
following equivalent statements are true:
\\
(a) The Young-Fenchel inequality
\[
\forall x\in X \colon \xi_r\of{x} \idif g\of{x} \leq g^*\of{\xi; r}.
\]
(b)
\[
\forall x\in X \colon \xi_r\of{x} \leq g\of{x} \isum g^*\of{\xi; r}.
\]
(c)
\[
\forall x\in X \colon \xi_r\of{x} \idif g^*\of{\xi; r} \leq g\of{x}.
\]
\end{proposition}

{\sc Proof.} The first equation is immediate from the definition on the conjugate. For
(b) and (c), recall $r \idif s \leq t$ is equivalent to $r \leq s \isum t$ (see theorem
\ref{ThmInfRes}, (a)). This gives the equivalence of (a) with (b) and (c), respectively.
\pend

\medskip Many of the known rules for the manipulation of conjugates apply also to $g^*$ at improper
elements. There are, however, some differences. We shall indicate one of them, a rule for
conjugates of an infimal convolution which is defined for $f, g \colon X \to \R^\triup$
(see \cite{Moreau63a}) as
\[
\of{f \square g}\of{x} = \inf\cb{f\of{x_1} \isum g\of{x_2} \colon x_1 + x_2 = x}.
\]

\begin{theorem}
\label{ThmInfConv} Let $f, g \colon X \to \R^\triup$ and $\xi \in X^\triup$, $r \in \R$.
Then
\begin{equation}
\label{EqInfConvGeneral} \of{f \square g}^*\of{\xi, r} =
   \sup\cb{ f^*\of{\xi, r_1} \ssum g^*\of{\xi, r_2} \colon r_1+r_2=r}.
\end{equation}
If, in particular, $\xi = x^* \in X^*$ then
\begin{equation}
\label{EqInfConvProper} \of{f \square g}^*\of{x^*, r} =
   f^*\of{x^*, 0} \ssum g^*\of{x^*, 0} \ssum \of{-r}.
\end{equation}
\end{theorem}

{\sc Proof.} From the first formula in (d) of proposition \ref{ProDiffCalc1} we obtain
\[
\of{f\square g}^*\of{\xi, r} = \sup\limits_{x,y \in X}\of{ \xi_r(x+y)\idif \of{f(x)\isum
g(y)}}.
\]
From \eqref{EqAffinSupConv1} and \eqref{EqConstPlusInf}, \eqref{EqProDualDiff1} we may
conclude
\begin{align*}
\of{f\square g}^*\of{\xi, r} & = \sup\limits_{x,y\in X}\of{\sup\limits_{r_1+r_2=r} \of{\xi_{r_1}(x)\ssum \xi_{r_2}(y)}\ssum (-1) \of{f(x)\isum g(y)} }\\
    & = \sup\limits_{\substack{x,y\in X\\r_1+r_2=x_0}}\Big(\big(\xi_{r_1}(x)\ssum \xi_{r_2}(y)\big)\ssum (-1) \big(f(x)\isum g(y)\big) \Big)
\end{align*}
Using proposition \ref{PropMinusSum} we get
\begin{align*}
\of{f\square g}^*\of{\xi, r}
    = & \sup\limits_{\substack{x,y\in X\\r_1+r_2=r}} \Big( \big(\xi_{r_1}(x)\ssum \xi_{r_2}(y)\big)\ssum \big((-1)f(x)\ssum (-1)g(y)\big) \Big)\\
    = & \sup\limits_{\substack{x,y\in X,\\ r_1+r_2=r}}\Big(\big(\xi_{r_1}(x)\ssum (-1)f(x)\big)\ssum \big(\xi_{r_2}(y)\ssum (-1)g(y)\big) \Big).
\end{align*}
Again \eqref{EqProDualDiff1} yields
\[
\of{f\square g}^*\of{\xi, r} =
    \sup\limits_{\substack{x, y \in X,\\ r_1+r_2=r}}
    \of{\of{\xi_{r_1}(x)\idif f(x)}\ssum \of{\xi_{r_2}(y)\idif g(y)}}.
\]
Taking the supremum with respect to $x$ while applying the second part of
\eqref{EqConstPlusInf}, and then doing the same with the supremum with respect to $y$ we
arrive at
\begin{align*}
\of{f\square g}^*\of{\xi, r}
    = & \sup\limits_{\substack{y\in X,\\ r_1+r_2=r}}\of{f^*\of{\xi, r_1}\ssum \of{\xi_{r_2}(y)\idif g(y)}}\\
    = & \sup\limits_{r_1+r_2=r}\of{f^*\of{\xi, r_1}\ssum g^*\of{\xi, r_2} }.
\end{align*}
This completes the proof. \pend

\begin{remark}
Taking $r = 0$ in \eqref{EqInfConvProper} and observing that $g^*\of{x^*, 0} =
g^*\of{x^*}$ is the classical Legendre-Fenchel conjugate we arrive at the correct version
of theorem 2.3.1, (ix) in \cite{Zalinescu02}. See already Moreau's paper \cite{Moreau66},
paragraph 6.h.
\end{remark}

\begin{remark}
Surprisingly, the right hand side of \eqref{EqInfConvGeneral} turns out to be the
supremal convolution of the conjugates of $f$ and $g$ -- with respect to the real
variable $r$. One may observe once again that conjugation via the formula given in
definition \ref{DefConjugate} changes the image space: $\ssum$ has to be used on the
right hand side of \eqref{EqInfConvGeneral} whereas $\isum$ appears on the left hand
side.
\end{remark}

\begin{definition}
\label{DefBiconjugate} The Legendre-Fenchel biconjugate of $g \colon X \to \R^\triup$ is
the function $g^{**} \colon X \to \R^\triup$ which is given by
\[
g^{**}\of{x} =
 \sup_{\of{\xi, r} \in X^\trido \times \R} \cb{\xi_r\of{x} \idif g^*\of{\xi, r}}
\]
for $x\in X$.
\end{definition}

\begin{theorem}\label{ThmBiconjugationScalar}
Let $g \colon X \to \R^\triup$ be closed convex. Then $g = g^{**}$.
\end{theorem}
{\sc Proof.} Note that for $x^* \in X^*$ and $r \in \R$ we have
\begin{equation}\label{BiconjFormulaImprop}
 \hat x_r^*\of{x} \idif g^*\of{\hat x^*, r} =
            \left\{
              \begin{array}{ccc}
                +\infty & : & x\not\in \dom \hat x_r^* \; \hbox{and} \;
                    \dom g\subseteq \dom \hat x_r^* \\
                -\infty & : & \hbox{otherwise}
              \end{array}
            \right.
\end{equation}
and
\begin{align}\label{BiconjFormulaProp}
 x_r^*\of{x} \idif g^*\of{x^*, r}=x^*\of{x} \idif g^*\of{x^*, 0}.
\end{align}

If $\epi g=\emptyset$, then $g^{**}(x)=g(x)=+\infty$ for all $x\in X$. If $\epi g \neq
\emptyset$, then $g$ is the pointwise supremum of its affine minorants. Especially, if
$g$ is improper, then it is the pointwise supremum of its improper affine minorants, thus
by \eqref{BiconjFormulaImprop} $g^{**}(x)=g\of{x}$ holds for all $x\in X$. If $g$ is
proper, then the well-known biconjugation theorem (see e.g. \cite[theorem
2.3.3]{Zalinescu02}) combined with \eqref{BiconjFormulaProp} delivers the desired result.
\pend

\begin{remark}
The closed convex hull of $g \colon X \to \R^\triup$ is defined by
\begin{align*}
  \forall x\in X \colon (\cl \co g)(x) = \inf \cb{t \in \R \colon \of{x,t} \in \cl\co\of{\epi g}}.
\end{align*}
Since the conjugate of an arbitrary function coincides with the conjugate of its closed
convex hull, the biconjugate of the function yields precisely the closed convex hull.
\end{remark}

\medskip The well-known relationship between the subdifferential and the Fenchel conjugate can be
extended to the improper case as follows.

\begin{proposition}\label{PropImpSubdiffConj}
Let $g \colon X \to \R^\triup$ be a convex function and $x_0 \in\dom g$. Then, the
following statements are equivalent for $\xi = \hat x^* \in \hat X^*$, $r=x^*(x_0)$:
\\
(a) $\xi \in \partial^{ex}g\of{x_0}$,
\\
(b) $\dom \xi_{r} \supseteq \dom g$,
\\
(c) $\forall x \in X$: $\of{\xi_r\of{x} \idif g\of{x}} \isum g\of{x_0} \leq \xi_r\of{0} =
-\infty$.
\\
(d) $g^*\of{\xi,r} = -\infty$.
\end{proposition}

{\sc Proof.} It holds $\xi_{x^*(x_0)}(x)=\xi(x-x_0)$ for all $x\in X$. Moreover,
$\xi_{x^*(x_0)}(x)=-\infty$ if $x^*(x)\leq x^*(x_0)$ and $\xi_{x^*(x_0)}(x)=+\infty$
otherwise. (a) $\Leftrightarrow$ (b) can be checked directly. (a) $\Leftrightarrow$ (d)
is formula \eqref{EqImpropConjug}. (d) $\Rightarrow$ (c) is clear from the definition of
$g^*$, while (c) and "not (b)" produce a contradiction. \pend

\begin{proposition}\label{PropSubdiffConj}
Let $g \colon X \to \R^\triup$ be a convex function and $x_0\in X$. If $x_0 \not\in \dom
g$, then $\partial^{ex}g\of{x_0} = \cb{-\infty}$. If $x_0 \in\dom g$, then
\begin{align*}
\partial^{ex}g\of{x_0} =
  \cb{\xi \in X^\triup \colon g^*\of{\xi, x^*\of{x_0}} \isum g\of{x_0} \leq \xi\of{0}}.
\end{align*}
\end{proposition}

{\sc Proof.} Assume $x_0 \in \dom g$. First, if $\xi=\hat x^* \in \hat X^*$ then
$\xi\of{0} = -\infty$. If $\xi \in \partial^{ex}g\of{x_0}$ then proposition
\ref{PropImpSubdiffConj}, (d) gives $g^*\of{\xi, x^*\of{x_0}} = -\infty$, hence
$g^*\of{\xi, x^*(x_0)} \isum g\of{x_0} \leq \xi\of{0}$. If $\xi$ satisfies the latter
inequality, then, by definition of $g^*$, also the one in proposition
\ref{PropImpSubdiffConj}, (c), hence $\xi \in \partial^{ex}g\of{x_0}$. Secondly, let
$\xi=x^* \in X^*$. Then $g$ is proper (otherwise $\partial g\of{x_0} =
\partial^{ex}g\of{x_0} \cap X^* = \emptyset$ and $g^*\of{\xi, x^*(x_0)} \equiv +\infty$),
and the result is well-known since $g^*\of{\xi, x^*\of{x_0}} = g^*\of{x^*}- x^*\of{x_0}$
in this case. \pend

\medskip We close this subsection by noting that the theory in sections 5.1, 5.2 has a symmetric
counterpart for functions mapping into $\R^\trido$ which requires modified definitions.
For example, Moreau \cite[p. 10]{Moreau63a}, already introduced the sup-convolution for
$\R^\trido$-valued functions.

\subsection{Set-valued convex and concave functions}

\subsubsection{Image spaces}

Let $Z$ be a topological linear space and $\mathcal{P}\of{Z}$ the collection of all
subsets of $Z$ including $\emptyset$. Let $C \subseteq Z$ be a convex cone including $0
\in Z$. We shall write $z_1 \leq_C z_2$ for $z_2 - z_1 \in C$, and this defines a
reflexive, transitive relation. The relation $\leq_C$ can be extended in two different
ways to $\mathcal{P}\of{Z}$, see for example \cite{Hamel05Habil}, \cite{Hamel09} for
details and references. A basic idea in these references is to use the equivalence
classes in $\mathcal{P}\of{Z}$ of the extension of $\leq_C$ to construct appropriate
image spaces for set-valued functions. It turned out (see \cite{Hamel09}) that the
following sets are appropriate choices as image spaces for set-valued closed convex
(concave) functions:
\begin{align*}
\mathcal{Q}^t_C\of{Z} & = \cb{A \subseteq Z \colon A = \cl\co\of{A + C}} \quad \mbox{and} \\
\mathcal{Q}^C_t\of{Z} & = \cb{A \subseteq Z \colon A = \cl\co\of{A -C}}.
\end{align*}
We redefine the addition for elements $A, B$ of $\mathcal{Q}^t_C\of{Z}$ and $\mathcal
{Q}^C_t\of{Z}$ and the multiplication with $0 \in \R$ by
\begin{equation}
\label{ClosedSum}
 A \oplus B = \cl\of{A+B}
\end{equation}
and $0 \cdot A = \cl C$ in $\mathcal{Q}^t_C\of{Z}$ and $0 \cdot A = -\cl C$ in
$\mathcal{Q}^C_t\of{Z}$, respectively. Then, $\of{\mathcal{Q}^t_C\of{Z},\oplus, \cdot,
\supseteq}$ and $\of{\mathcal{Q}^C_t\of{Z}, \oplus, \cdot, \subseteq}$ are partially
ordered conlinear spaces. Again, the multiplication with $-1$, transforming
$\mathcal{Q}^t_C\of{Z}$ into $\mathcal{Q}^C_t\of{Z}$ and vice versa, is a duality in the
sense of \cite{Singer97}. We shall abbreviate $\mathcal{Q}^{\triup} =
\of{\mathcal{Q}^t_C\of{Z},\oplus, \cdot, \supseteq}$ and $\mathcal{Q}^{\trido} =
\of{\mathcal{Q}^C_t\of{Z},\oplus, \cdot, \subseteq}$. The advantage of using  these image
spaces compared to other, more common approaches in vector and set optimization is that
they are partially ordered lattices with formulas for $\inf$ and $\sup$ as given below.
Moreover, $\inf$ and $\sup$ in these spaces are not "utopia elements", but they are
strongly related to known extremality concepts based on minimal/maximal points with
respect to $\leq_C$ or set relations. Compare \cite{Hamel05Habil}, \cite{Hamel10} and,
most notably, the discussion in \cite{HeydeLoehne09R}.

\begin{proposition}
(a) Let $\mathcal A \subseteq \mathcal{Q}^{\triup}$ and $B \in \mathcal{Q}^{\triup}$.
Then
\[
 \inf \mathcal A = \cl\co\bigcup\limits_{A\in \mathcal A}A, \quad
 \sup \mathcal A = \bigcap\limits_{A \in \mathcal A} A
\]
and
\[
 \inf \mathcal A \oplus B = \inf\of{\mathcal A \oplus \cb{B}}, \quad
 \sup \mathcal A \oplus B \subseteq \sup\of{\mathcal A \oplus \cb{B}}.
\]
(b) Let $\mathcal A \subseteq \mathcal \mathcal{Q}^{\trido}$ and $B\in \mathcal
\mathcal{Q}^{\trido}$. Then
\[
 \sup \mathcal A = \cl\co\bigcup\limits_{A \in \mathcal A} A, \quad
 \inf \mathcal A = \bigcap\limits_{A \in \mathcal A} A
\]
and
\[
 \sup \mathcal A \oplus B = \sup\of{\mathcal A \oplus \cb{B}}, \quad
 \inf \mathcal A \oplus B \subseteq \inf\of{\mathcal A \oplus \cb{B}}.
\]
The addition for sets of sets is defined as $\mathcal A \oplus \cb{B} = \cb{A \oplus B
\colon A \in \mathcal A}$ if $\mathcal A$ is non-empty and $\mathcal A \oplus B =
\emptyset$ otherwise.
\end{proposition}

{\sc Proof.} Exemplarily, we add the proof of infimum-additivity in $Q^\triup$. Indeed,
for $\mathcal A \subseteq \mathcal{Q}^{\triup}$ and $B \in \mathcal{Q}^{\triup}$ we have
\[
\inf \mathcal A \oplus B = \of{\cl\co\bigcup\limits_{A\in \mathcal A}A} \oplus B =
    \of{\co\bigcup\limits_{A\in \mathcal A}A} \oplus B =
    \cl\co\bigcup\limits_{A\in \mathcal A}\of{A \oplus B} = \inf\of{\mathcal A \oplus
    \cb{B}}.
\]
Since the roles of $\inf$ and $\sup$ are exchanged in $Q^\trido$, the supremum-additivity
in $Q^\trido$ follows directly. The other results are essentially a consequence of the
definitions of $\inf$, $\sup$ and $\oplus$. \pend

\medskip In view of the theorems \ref{ThmInfRes}, \ref{ThmSupRes}, the previous proposition
tells us that $\mathcal{Q}^{\triup}$ and $\mathcal{Q}^{\trido}$ are order complete
residuated lattices. Note that $\emptyset$ is the greatest element in
$\mathcal{Q}^{\triup}$ and the least in $\mathcal{Q}^{\trido}$. In both cases, it
dominates the addition which is in complete analogy with $+\infty$ dominating the
inf-addition in $\R^\triup$ and $-\infty$ the sup-addition in $\R^\trido$. The
residuation can be used to define a difference for sets.

\begin{definition}\label{DefSetDifference}
Let $A, B \in Q^\triup$. The set
\begin{align}
  A \idif B = \inf\cb{M\in Q^\triup \colon A \supseteq B\oplus M} = \cb{z\in Z \colon B+\cb{z} \subseteq A}
\end{align}
is called the inf-difference of $A$ and $B$. Likewise, for $A, B\in Q^\trido$ the set
\begin{align}
A \sdif B = \sup\cb{M\in Q^\trido \colon B\oplus M \subseteq A} =
    \cb{z\in Z \colon B+\cb{z}\subseteq A}
\end{align}
is called the sup-difference of $A$ and $B$.
\end{definition}

The reader may wonder why we use the same expression for the difference in
$\mathcal{Q}^{\triup}$ and $\mathcal{Q}^{\trido}$ which is not the case in $\R^\triup$
and $\R^\trido$. Of course, the reason is that we use different order relations, namely
$\supseteq$ in $\mathcal{Q}^{\triup}$ and $\subseteq$ in $\mathcal{Q}^{\trido}$ and
therefore, the infimum in $\mathcal{Q}^{\triup}$ is a union as it is the supremum in
$\mathcal{Q}^{\trido}$. We had the same effect if we would use $\leq$ in $\R^\triup$ and
$\geq$ in $\R^\trido$. However, note that $A \idif B \in \mathcal{Q}^{\triup}$ for $A, B
\in \mathcal{Q}^{\triup}$ while $A \sdif B \in \mathcal{Q}^{\trido}$ for $A, B \in
\mathcal{Q}^{\trido}$, so $A \idif B$ and $A \sdif B$ look very different in general.

\medskip From now on, let $X$ and $Z$ be separated, locally convex spaces with topological
duals $X^*$ and $Z^*$, respectively. The (negative) dual cone of $C$ is the set
\[
C^- = \cb{z^* \in Z^* \colon \forall z \in C \colon z^*\of{z} \leq 0}.
\]
It is well-known (and a consequence of a separation argument) that sets $A \in
\mathcal{Q}^{\triup}$, $B \in \mathcal{Q}^{\trido}$ can be described dually as
\begin{align*}
A & = \bigcap\limits_{z^*\in C^-\bs\cb{0}}
    \cb{z\in Z \colon \inf\limits_{a\in A} -z^*\of{a} \leq -z^*\of{z}}, \\
B & = \bigcap\limits_{z^*\in C^-\bs\cb{0}}
    \cb{z \in Z \colon  -z^*\of{z} \leq \sup\limits_{b\in B} -z^*\of{b}},
\end{align*}
respectively. The above representation of elements in $\mathcal{Q}^\triup$ and
$\mathcal{Q}^\trido$ can be used to characterize the set differences from definition
\ref{DefSetDifference} in terms of support functions. For $D \subseteq Z$, define the
extended real-valued functions $\sigma^\triup_D \colon Z^* \to \R^\triup$,
$\sigma^\trido_D \colon Z^* \to \R^\trido$ by
\[
\sigma^\triup_D\of{z^*} = \inf_{z \in D}-z^*\of{z}, \quad
    \sigma^\trido_D\of{z^*} = \sup_{z \in D}-z^*\of{z}.
\]

\begin{proposition}
\label{PropScalarSetDiff} (a) For all $A, B\in \mathcal Q^\triup$,
\[
A \idif B = \bigcap\limits_{z^*\in C^-\bs\cb{0}}
    \cb{z\in Z \colon \sigma^\triup_A\of{z^*} \idif \sigma^\triup_B\of{z^*} \leq -z^*\of{z}}.
\]
In particular, if $A = \cb{z \in Z \colon \sigma^\triup_A\of{z^*} \leq -z^*(z)}$ for
$z^*\in C^-\bs\cb{0}$, then
\begin{align*}
A \idif B = \cb{z \in Z \colon \sigma^\triup_A\of{z^*} \idif \sigma^\triup_B\of{z^*} \leq
-z^*\of{z}}.
\end{align*}
(b) For all $A, B\in \mathcal Q^\trido$,
\[
A \sdif B = \bigcap\limits_{z^*\in C^-\bs\cb{0}}
    \cb{z \in Z \colon -z^*\of{z}\leq \sigma^\trido_A\of{z^*} \sdif
    \sigma^\trido_B\of{z^*}}.
\]
In particular, if $A=\cb{z \in Z\colon -z^*(z) \leq \sigma^\trido_A\of{z^*}}$ for $z^*\in
C^-\bs\cb{0}$, then
\begin{align*}
A \sdif B =
    \cb{z\in Z\colon -z^*(z) \leq  \sigma^\trido_A\of{z^*} \sdif \sigma^\trido_B\of{z^*}}.
\end{align*}
\end{proposition}

{\sc Proof.} (a) By definition of $A \idif B$, $\sigma^\triup_A\of{z^*} \leq
\sigma^\triup_{B+z}\of{z^*}$ whenever  $z \in A \idif B$, and
$\sigma^\triup_{B+z}\of{z^*} = \sigma^\triup_B\of{z^*} \isum z^*\of{-z}$ for all $z^*\in
C^-\bs\cb{0}$. Thus, $\sigma^\triup_A\of{z^*} \idif \sigma^\triup_B\of{z^*} \leq
-z^*(z)$.

On the other hand, take $z\in Z$ such that $\sigma^\triup_A\of{z^*} \idif
\sigma^\triup_B\of{z^*} \leq -z^*\of{z}$ for all $z^*\in C^-\bs\cb{0}$ and assume $z
\notin A\idif B$. Then there is $b \in B$ such that $z+b \notin A$. A separation argument
produces $z^*\in C^-\bs\cb{0}$ and $t \in \R$ such that $-z^*\of{b+z}< t \leq
\sigma^\triup_A\of{z^*}$. Since $\sigma^\triup_B\of{z^*} \leq -z^*\of{b}$ we may conclude
$-z^*\of{z} < \sigma^\triup_A\of{z^*} \idif \sigma^\triup_B\of{z^*}$, a contradiction.

If, additionally, $A = \cb{z\in Z\colon \sigma^\triup_A\of{z^*} \leq -z^*\of{z}}$ with
$z^* \in C^-\bs\cb{0}$ and $\sigma^\triup_A\of{z^*} \idif \sigma^\triup_B\of{z^*} \leq
-z^*\of{z}$ for $z\in Z$, then $\sigma^\triup_A\of{z^*} \leq -z^*\of{b+z}$ for all $b \in
B$, thus $B+z\subseteq A$.

(b) is proven with parallel arguments. \pend

\medskip From proposition \ref{PropScalarSetDiff} we immediately obtain for all
$A, B \subseteq \mathcal Q^\triup$ and $z^*\in C^-\bs\cb{0}$
\[
\sigma^\triup_{A \idif B}\of{z^*} \geq \sigma^\triup_A\of{z^*} \idif
\sigma^\triup_B\of{z^*}
\]
with equality if $A = \cb{z \in Z \colon \sigma^\triup_A\of{z^*} \leq -z^*(z)}$. A
parallel relationship holds for elements of $\mathcal Q^\trido$ and $\sdif$.

\begin{remark} For $z^* \in Z^*$ we set $H\of{z^*} = \cb{z \in Z \colon z^*\of{z} \leq
0}$. The additional assumption in proposition \ref{PropScalarSetDiff}, (a) is equivalent
to $A \oplus H\of{z^*} = A$ since
\[
A \oplus H\of{z^*} = \cb{z\in Z \colon \vp^\triup_{z^*}\of{A}\leq -z^*(z)}.
\]
\end{remark}

\subsubsection{Set-valued functions and their scalar representation}

The graph of a function $f \colon X \to \mathcal P\of{Z}$ is the set
\[
\gr f = \cb{\of{x, z} \in X \times Z \colon z \in f\of{x}},
\]
and the effective domain of $f$ is $\dom f = \cb{x \in X \colon f\of{x} \neq \emptyset}$.

\begin{definition}
A function $g \colon X \to \mathcal{Q}^{\triup}$ is called convex (positively homogenous,
closed) iff $\gr g \subseteq X \times Z$ is convex (a cone, closed). A positively
homogenous convex function into $\mathcal{Q}^{\triup}$ is called sublinear. A function $h
\colon X \to \mathcal{Q}^{\trido}$ is called concave (positively homogenous, closed) iff
$\gr h \subseteq X \times Z$ is a convex (a cone, closed). A positively homogenous
concave function into $\mathcal{Q}^{\trido}$ is called superlinear.
\end{definition}

The collection of all convex functions into $\mathcal{Q}^{\triup}$
($\mathcal{Q}^{\trido}$) is a conlinear space under point-wise addition and
multiplication with non-negative reals. As in the scalar case, it does not make much
sense to consider convex functions into $\mathcal{Q}^{\trido}$ or concave into
$\mathcal{Q}^{\triup}$. A function $g \colon X \to \mathcal{Q}^{\triup}$ is convex if and
only if $-g \colon X \to \mathcal{Q}^{\trido}$ is concave.

The next goal is to represent $\mathcal{Q}^{\triup}$-valued functions by families of
$\R^\triup$-valued ones and, likewise, $\mathcal{Q}^{\trido}$-valued functions by
families of $\R^\trido$-valued ones. Let $z^*\in C^-$ and $g \colon X \to
\mathcal{Q}^\triup$, $h \colon X \to \mathcal{Q}^\trido$. Define $\vp_{g, z^*}^\triup
\colon X \to \R^\triup$ and $\vp_{h, z^*}^\trido \colon X \to \R^\trido$, respectively,
by
\begin{equation}
\label{EqScalarFamilies} \vp_{g, z^*}^\triup\of{x} =
    \inf_{z\in g\of{x}} -z^*\of{z} \quad \hbox{and} \quad
    \vp_{h, z^*}^\trido\of{x} = \sup_{z\in h\of{x}} - z^*\of{z}
\end{equation}
for $x \in X$. Note that $g$ is convex if and only if $\vp_{g, z^*}^\triup$ is convex for
all $z^*\in C^-$, and $h$ is concave if and only if $\vp_{h, z^*}^\trido$ is concave for
all $z^*\in C^-$, see \cite[lemma 3.2.3.]{Schrage09Diss}. If $z^*=0$ then
\[
\vp_{g, z^*}^\triup\of{x} =
    \left\{
    \begin{array}{ccc}
    0 & : & x\in \dom g \\
    +\infty & : & \hbox{otherwise.}
    \end{array}
    \right.
\]
Moreover, $-\vp_{-g, z^*}^\trido\of{x} = \vp_{g, z^*}^\triup\of{x}$ for all $x \in X$.
From the dual description of elements of $\mathcal{Q}^\triup$, $\mathcal{Q}^\trido$ the
following formulas are immediate for functions $g \colon X \to \mathcal{Q}^\triup$, $h
\colon X \to \mathcal{Q}^\trido$:
\begin{align} \forall x \in X \colon & g\of{x} =
    \bigcap\limits_{z^*\in C^-\bs\cb{0}}
    \cb{z \in Z \colon \vp_{g, z^*}^{\triup}\of{x} \leq -z^*\of{z}}\\
\forall x \in X \colon & h\of{x} =
    \bigcap\limits_{z^*\in C^-\bs\cb{0}}
    \cb{z \in Z \colon -z^*\of{z} \leq \vp_{h, z^*}^{\trido}\of{x}}.
\end{align}
These formulas tell us that $\mathcal{Q}^\triup$- and $\mathcal{Q}^\trido$-valued
functions can be represented by families of extended real-valued functions. In general,
the scalarizations may behave "very badly".

\begin{example}\label{ScalarizationNotClosed}
Let the set $\R^2$ be ordered by the cone $C=\R^2_+$, $z^*=(0,-1)$ and $g:\R\to\R^2$ be
defined as $g(x)=\cb{(\frac{1}{x},0)}+C$ if $x>0$, and $g(x)=\emptyset$ otherwise. Then,
$\vp^\triup_{g,z^*}(0)=+\infty$, while $\vp^\triup_{g,z^*}(x)=0$ holds for all $x>0$ and
thus $\cl \vp^\triup_{g,z^*}(0)=0$.
\end{example}

The above example shows that the scalarizations $\vp_{g, z^*}^\triup$ of a closed
function need not be closed. However, one can restrict the scalarizations to closed ones
as already shown in \cite[proposition 3.3.5]{Schrage09Diss}.

\begin{proposition}
\label{PropClosedScalarization} Let $g \colon X \to \mathcal{Q}^\triup$ and $h \colon X
\to \mathcal{Q}^\trido$ be a closed convex and a closed concave function, respectively.
Then
\begin{align}
\label{EqScalar1} \forall x \in X \colon g\of{x} & = \bigcap\limits_{z^*\in C^-\bs\cb{0}}
    \cb{z\in Z \colon \of{\cl \vp_{g, z^*}^\triup}\of{x} \leq -z^*\of{z}} \\
\label{EqScalar2} \forall x \in X \colon h\of{x} & = \bigcap\limits_{z^*\in C^-\bs\cb{0}}
    \cb{z \in Z \colon -z^*\of{z} \leq \of{\cl \vp_{h, z^*}^{\trido}}\of{x}}.
\end{align}
\end{proposition}

{\sc Proof.} A function $h \colon X \to Q^\trido$ is closed and concave if and only if
$-h \colon X \to Q^\triup$ is closed and convex. Thus, it is sufficient to prove the
statement for convex functions. Let $g \colon X\to Q^\triup$ be closed and convex. If
$\gr g = \emptyset$, then there is nothing to prove. Let $\of{x_0,z_0} \notin \gr g$.
Then by a separation argument in $X \times Z$ there exists $\of{x^*, z^*} \in \of{X^*
\times Z^*} \bs\cb{\of{0,0}}$ and $t \in \R$ such that
\[
-x^*(x_0)-z^*(z_0)< t < \inf\limits_{(x,z) \in \gr g}\of{-x^*(x)-z^*(z)}.
\]
Obviously, $z^*\in C^-$ and
\[
\forall x\in X \colon \quad x^*(x)-\of{x^*(x_0)+z^*(z_0)}<x^*(x)+t < \inf\limits_{z\in
g(x)}\of{-z^*(z)}=\vp^\triup_{g,z^*}(x).
\]
Thus, $x\mapsto x^*(x)+t$ is an affine minorant of $\cl \vp^\triup_{g,z^*}$ and
\[
z_0\notin \cb{z\in Z:\; \vp^\triup_{g,z^*}(x_0)\leq -z^*(z)}.
\]
This proves $\supseteq$ in \eqref{EqScalar1}, and the converse inclusion is trivial.
\pend

\medskip A function $g \colon X \to \mathcal{Q}^\triup$ is called proper  ($C$-proper) if
$\dom f \neq \emptyset$ and $f\of{x} \neq Z$ ($f\of{x} \neq f\of{x} - C$) holds for all $
\in X$. A function $g \colon X \to \mathcal{Q}^\triup$ is proper ($C$-proper) if and only
if there exists at least one $z^* \in C^-\bs\cb{0}$ ($z^*\in C^-\bs -C^-$) such that
$\vp_{g, z^*}^\triup \colon X \to \R^\triup$ is proper. However, even if $g$ is
$C$-proper, not all scalarizations are proper in general.

\begin{example}\label{ImproperScalarization1}
Let $Z = \R^2$, $C = \cl\cone\cb{c}$ with $c = \of{0,1}^T$. Take $z^*_0=(0,-1)^T \in
\R^2$ and define a function $g \colon \R \to \mathcal{Q}^\triup$ by
\begin{equation*}
g(x)=\left\{
     \begin{array}{ccc}
     H(z^*_0) = \cb{z \in \R^2 \colon \of{z^*_0}^Tz \leq 0} & : & x > 0\\
     C & : & x=0 \\
     \emptyset & : & x < 0
     \end{array}
     \right.
\end{equation*}
Then $g$ is convex $C$-proper, and $\vp_{g, z^*}^\triup\colon X \to \R^\triup$ is proper
if and only if $z^*$ is collinear with $z^*_0$. Moreover, the function $g$ is not
completely characterized by its proper scalarizations since $g(0)=C\subsetneq H(z^*_0)$
holds.
\end{example}

The previous example shows that while analyzing set-valued functions, improper scalar
functions appear naturally. Providing a calculus on the space of improper functions
allows a unified approach to the theory of set-valued functions via scalarization.

\subsubsection{Conaffine proper and improper functions}

In the remaining two sections we focus on the convex case mentioning that the
corresponding constructions for the concave one are easily obtained.

\begin{definition}\label{DefConlinearFunction}
Let $x^* \in X^*$, $z^* \in C^-$. The function $S^\triup_{\of{x^*,z^*}} \colon X \to
\mathcal{Q}^\triup$ given by
\[
S^\triup_{\of{x^*,z^*}}\of{x} = \cb{z \in Z \colon x^*\of{x} + z^*\of{z} \leq 0}
\]
for $x\in X$ is called an lower conlinear function.
\end{definition}

Obviously, $S^\trido_{\of{x^*,z^*}}\of{x} = -S^\triup_{\of{x^*,z^*}}\of{x}$ for $x\in X$
is the corresponding upper conlinear function to be used in the concave case. For each
$z^* \neq 0$, the functions $S^\triup_{\of{x^*,z^*}}$, $S^\trido_{\of{x^*,z^*}}$ are
"finite-valued" in the sense that they attain neither the value $Z$ nor $\emptyset$.

It is easy to find situations in which $S^\triup_{\of{x^*,z^*}}$ is $C$-proper, but its
scalarization with $z^*_0$ is improper (in which case $z^*_0$ is not collinear with
$z^*$). Of course, the scalarization with $z^*$ is linear if $z^* \neq 0$. If $z^* = 0$,
then the scalarization of a lower conlinear function is $\hat x^*$.

For each $z^* \in Z^*\bs\cb{0}$, there is a one-to-one correspondence between functions
$x^* \in X^*$ and $S^\triup_{(x^*,z^*)} \colon X \to Q^\triup$. The situation is
different for $z^* = 0$, see remark \ref{RemZeroImpropConaffine} below.

In the same way as we extended the definition of affine functions from real-valued to
extended real-valued ones we extend the definition of conlinear functions.

\begin{definition}
\label{DefImpropConaffFunc} Let $\xi \in X^\triup$ and $r \in \R$. We define the
functions $S^\triup_{(\xi, r, z^*)} \colon X \to Q^\triup$ by
\[
S^\triup_{(\xi, r, z^*)}(x) = \cb{z \in Z \colon \xi_r(x) \leq -z^*(z)}
\]
for $x \in X$. Each such function is called a (closed) conaffine function. If $r = 0$ it
is called a (closed) conlinear function.
\end{definition}

Of course, a $Q^\triup$-valued conaffine function is proper if and only if $z^* \neq 0$
and $\xi \in X^*$.

\begin{remark}
\label{RemZeroImpropConaffine} Let $z^* \in Z^*$, $\xi = \hat x^* \in \hat X^*$ and $r
\in \R$. Then
\begin{align*}
\forall x \in X \colon
    S^\triup_{(\hat x^*, r, z^*)}(x) & =
    \cb{z\in Z:\; \hat x_r^*(x) \leq -z^*(z)} \\
    & =
    \left\{
    \begin{array}{lll}
    Z & : & x_r^*\of{x} \leq 0 \; (\Leftrightarrow \; x \in \dom \hat x^*_r) \\
    \emptyset & : & x_r^*\of{x} > 0 \; (\Leftrightarrow \; x \not\in \dom \hat x^*_r)
    \end{array}
    \right. \\
    & = S^\triup_{(x^*, r, 0)}(x).
\end{align*}
This means, one can replace $\hat x^*$ by $x^*$ and $z^*$ by $0 \in Z^*$, and in this way
the case $z^* = 0$ "includes" all improper cases. In particular,
\[
\forall x \in X \colon
    S^\triup_{(x^*,0)}\of{x} = \cb{z \in Z \colon x^*\of{x} \leq 0} =
    \left\{
    \begin{array}{lll}
    Z & : & x^*\of{x} \leq 0 \; (\Leftrightarrow \; x \in \dom \hat x^*) \\
    \emptyset & : & x^*\of{x} > 0 \; (\Leftrightarrow \; x \not\in \dom \hat x^*)
    \end{array}
    \right.
\]
\end{remark}

Consequently, we obtain a one-to-one correspondence between the set of improper affine
scalar functions $\hat x_r^* \colon X \to \R^\triup$ and the set of conaffine functions
$S^\triup_{(\hat x^*, r, 0)} \colon X \to Q^\triup$. Note that many of those scalar and
set-valued functions coincide since $\hat x^* \equiv t \hat x^*$ for $t>0$.

Finally, we turn to scalarizations of proper and improper conaffine functions.

\begin{proposition}
\label{PropScalarConaffine} Let $z^*, z^*_0 \in Z^*\bs\cb{0}$, $x^* \in X^*$, $r\in\R$.
Then
\begin{align}
\vp^\triup_{(S^\triup_{(x^*, r, z^*_0)},z^*)}\of{x} & =
    \left\{
    \begin{array}{ccc}
    -\infty & : & \forall t > 0 \colon z^* \neq tz^*_0 \\
    tx^*\of{x - x_0} & : & z^* = tz^*_0, \; t > 0
    \end{array}
    \right. \\
\vp^\triup_{S^\triup_{(\hat x^*, r, 0)}, z^*}\of{x} & =
    \vp^\triup_{S^\triup_{(x^*, r, 0)}, z^*}\of{x} =
    \hat x_r^*\of{x}.
\end{align}
\end{proposition}

{\sc Proof.} Obvious from the definition in \eqref{EqScalarFamilies}. \pend

\subsubsection{Legendre--Fenchel conjugates of set-valued functions}

Let $\xi \in X^\triup$, $r \in \R$, $z^* \in C^-\bs\cb{0}$. The function
$S^\triup_{\of{\xi, r, z^*}} \colon X \to Q^\triup$ is a closed conaffine minorant of $g
\colon X \to Q^\triup$ if and only if
\begin{align}\label{EqConaffMinorant}
\forall x\in X \colon S^\triup_{\of{\xi, r, z^*}}\of{x} \supseteq g\of{x}.
\end{align}
The following result runs parallel to theorem \ref{ThmImpropAffMin}. It should be clear
how to formulate the counterpart for concave functions.

\begin{theorem}
\label{ThmImpropAffMinSet} Let $g \colon X \to \mathcal Q^\triup$ be improper closed
convex. Then $g$ is the pointwise supremum of its improper closed conaffine minorants,
that is
\begin{equation}
\label{EqAffMinSet1} \forall x \in X \colon
    g\of{x} = \bigcap\cb{S^\triup_{\of{\hat x^*, r, z^*}}\of{x} \colon
    \of{\hat x^*, r, z^*} \in \hat X^* \times \R \times C^-\bs\cb{0} \colon
    \text{\eqref{EqConaffMinorant} is satisfied}}.
\end{equation}
In particular, if $g \not\equiv Z$, then there are $x^* \in X^*\bs\cb{0}$, $z^* \in
C^-\bs\cb{0}$ and $r \in \R$ such that
\begin{equation}
\label{EqAffMinSet2} \forall x \in X \colon S^\triup_{\of{\hat x^*, r, z^*}}\of{x}
\supseteq g\of{x}.
\end{equation}
\end{theorem}

{\sc Proof.} The theorem is trivial if $g \equiv Z$ and $g \equiv \emptyset$. Let us
assume that $g$ is different from these two functions. In this case, $\dom g \neq
\emptyset$ is a closed convex set and $g = Z$ on $\dom g$ (see proposition 5 in
\cite{Hamel09}). Hence $\dom g$ is the intersection of non-trivial closed half spaces
including it. Each such half space is generated by some $x^* \in X^*\bs\cb{0}$, $r \in
\R$, and it is easily seen with the help of remark \ref{RemZeroImpropConaffine} that the
corresponding $S^\triup_{\of{\hat x^*, r, z^*}}$ satisfies \eqref{EqAffMinSet2} with an
arbitrary $z^* \in C^-\bs\cb{0}$. \pend

\medskip The "proper/$C$-proper" version of this theorem is theorem 1 in \cite{Hamel09}:
$g \colon X \to Q^\triup$ is closed, convex and proper, or identically $Z$ or
$\emptyset$, if and only if $g$ is the pointwise supremum of its proper conaffine
minorants (that is, with $z^*\in C^-\bs\cb{0}$). Therefore, a closed convex function with
values in $\mathcal Q^\triup$ is the pointwise supremum of its conaffine closed
minorants. The next result is the set-valued counterpart of theorem \ref{ThmImpAffMin}.

\begin{theorem}
\label{ThmImpAffMinSet} Let $g \colon X \to \mathcal Q^\triup$ and $\xi \in X^\triup$,
$z^* \in C^-\bs\cb{0}$. The following statements are equivalent:
\\
(a) $\forall x \in X$: $S^\triup_{\of{\xi, r, z^*}}\of{x} \supseteq g\of{x}$, i.e.
$S^\triup_{\of{\xi, r, z^*}}$ is a closed conaffine minorant of $g$,
\\
(b) $\forall x \in X$: $\xi_r\of{x} \leq \vp^\triup_{g, z^*}\of{x}$, i.e. $\xi_r$ is an
affine minorant of $\vp^\triup_{g,z^*}$,
\\
(c) $\of{\vp^\triup_{g, z^*}}^*\of{\xi, r} \leq 0$,
\\
(d) $\sup \cb{S^\triup_{\of{\xi, r, z^*}}\of{x} \idif g\of{x} \colon x\in X} \supseteq
H\of{z^*}$.

Moreover, if $\xi = \hat x^*$ with $x^* \in X^*$, then (a) through (d) are equivalent to
\\
(e) $\dom \hat x_r^*\supseteq \dom g$,
\\
and the supremum in (d) is $Z$.
\end{theorem}

{\sc Proof.} (a)$\Rightarrow$(b): By assumption,
\begin{equation}
\label{EqThmImpAffMinSet1} \forall x \in X \colon
 \vp^\triup_{S^\triup_{\of{\xi, r, z^*}}, z^*}\of{x} \leq \vp^\triup_{g, z^*}\of{x}.
\end{equation}
If $z^* \neq 0$ then, by definition of $S^\triup_{\of{\xi, r, z^*}}$,
\[
\vp^\triup_{S^\triup_{\of{\xi, r, z^*}}, z^*}\of{x} =
    \inf\cb{-z^*\of{z} \colon \xi_r\of{x} \leq -z^*\of{z}} = \xi_r\of{x}
\]
for each $x \in X$. Now, (b) follows from \eqref{EqThmImpAffMinSet1}. If $z^* = 0$ then
$\vp^\triup_{g, 0}\of{x} = I_{\dom g}\of{x}$ and
\[
\vp^\triup_{S^\triup_{\of{\xi, r, 0}}, 0}\of{x} =
    \left\{
    \begin{array}{ccc}
    +\infty & : & \xi_r\of{x} > 0 \\
    0 & : & \xi_r\of{x} \leq 0
    \end{array}
    \right.
\]
With this, (b) is immediate from \eqref{EqThmImpAffMinSet1}.

(b)$\Rightarrow$(c): We have
\[
\of{\vp^\triup_{g, z^*}}^*\of{\xi, r} =
    \sup_{x \in X}\cb{\xi_r\of{x} \idif \vp^\triup_{g, z^*}\of{x}}
    \leq \sup_{x \in X}\cb{\xi_r\of{x} \idif \xi_r\of{x}} \leq 0
\]
where the first inequality follows with the help of proposition \ref{ProDiffCalc1} (b),
second formula, and the second inequality with the help of proposition \ref{ProDiffCalc1}
(a).

(c)$\Rightarrow$(d): If $z^* \neq 0$ we obtain from proposition \ref{PropScalarSetDiff}
(a) with $A = S^\triup_{\of{\xi, r, z^*}}\of{x}$, $\sigma^\triup_A\of{z^*} = \xi_r\of{x}$
and $B = g\of{x}$
\[
S^\triup_{\of{\xi, r, z^*}}\of{x} \idif g\of{x} =
    \cb{z \in Z \colon \xi_r\of{x} \idif \vp^\triup_{g, z^*}\of{x} \leq -z^*\of{z}}.
\]
Consequently,
\begin{multline*}
\bigcap_{x \in X} \sqb{S^\triup_{\of{\xi, r, z^*}}\of{x} \idif g\of{x}} =
    \cb{z \in Z \colon \sup_{x \in X} \sqb{\xi_r\of{x} \idif \vp^\triup_{g, z^*}\of{x}} \leq -z^*\of{z}}
\\ = \cb{z \in Z \colon \of{\vp^\triup_{g, z^*}}^*\of{\xi, r} \leq -z^*\of{z}}.
\end{multline*}
The desired implication follows since $H\of{z^*} = \cb{z \in Z \colon 0 \leq
-z^*\of{z}}$.

(d)$\Rightarrow$(a): By assumption, $H\of{z^*} \subseteq \cb{z \in Z \colon g\of{x} +
\cb{z} \subseteq S^\triup_{\of{\xi, r, z^*}}\of{x}}$ for all $x \in X$. Since $0 \in
H\of{z^*}$, this implies $g\of{x} \subseteq S^\triup_{\of{\xi, r, z^*}}\of{x}$ for all $x
\in X$ which is (a).

Finally, if $\xi = \hat x^*$, the equivalence of (a) and (e) follows from remark
\ref{RemZeroImpropConaffine}. \pend

\medskip The previous theorem motivates the following definition of set-valued conjugates for functions
$g \colon X \to Q^\triup$. In contrast to \cite{Hamel09}, we define $Q^\triup$-valued
conjugates via a supremum instead negative conjugates via an infimum. The definition
below is also slightly different from the one given in \cite{Schrage09Diss}, but uses the
same basic idea, namely the set difference.

\begin{definition}\label{DefSetValuedConjugate}
The Legendre-Fenchel conjugate $g^* \colon X^\triup \times \R \times C^- \to Q^\triup$ of
a function $g \colon X \to Q^\triup$ is the function $g^* \colon X^\triup \times \R
\times C^- \to Q^\triup$ defined by
\begin{align*}
g^*\of{\xi, r , z^*} =
    \bigcap\limits_{x\in X} \of{S^\triup_{(\xi, r, z^*)}(x) \idif g(x)}.
\end{align*}
\end{definition}

By remark \ref{RemZeroImpropConaffine} it holds
\begin{align}\label{EqProperDualVariables}
\forall (x^*, r,z^*) \in X^* \times \R \times C^-:\quad g^*(\hat x^*, r , z^*) =
g^*(x^*,r,0^*).
\end{align}

The definition of the conjugate together with theorem \ref{ThmImpAffMinSet} produces the
following scalarization formula for conjugates
\begin{align}
\label{EqScalarizationOfConjugate} g^*(\xi,r, z^*) = \cb{z\in Z \colon
\of{\vp^\triup_{g,z^*}}^*\of{\xi, r} \leq -z^*(z)}
\end{align}
for all $\xi \in X^\triup$, $r \in \R$, $z^* \in C^-$.

Most rules for manipulating conjugates carry over from the scalar case. In particular, if
$r = 0$, $\xi = x^* \in X^*$ and $z^*\in C^-\bs\cb{0}$, then all classic duality results
from the scalar theory can be proven for the set-valued case as well, compare
\cite{Hamel09}, \cite{Schrage09Diss}, \cite{Hamel10}.

We will close this section illustrating the previous statement using the biconjugation
theorem as an example. If $g \colon X \to Q^\triup$ is a function, its biconjugate is
defined to be
\begin{align}\label{BiconjugateSetValued}
  \forall x \in X \colon
 g^{**}(x) = \bigcap\limits_{\xi \in X^\triup, \, r \in \R,\, z^*\in C^-}
 \of{S^\triup_{(\xi, r, z^*)}\of{x} \idif g^*\of{\xi, r ,z^*}}.
\end{align}
The function $g^{**} \colon X \to Q^\triup$ maps indeed into $Q^\triup$ and is closed and
convex. By equation \eqref{EqProperDualVariables}, it holds
\begin{align}\label{BiconjugateSetValued3}
  \forall x \in X:\quad
 g^{**}(x) = \bigcap\limits_{\substack{(\xi,r) \in X^\triup \times \R,\\ z^*\in C^-\setminus\cb{0}}}
 \of{S^\triup_{(\xi, r, z^*)}\of{x} \idif g^*\of{\xi, r ,z^*}}.
\end{align}

\begin{theorem}\label{ThmBiconjugateSetValued}
Let $g \colon X\to Q^\triup$ be closed and convex. Then $g^{**} = g$.
\end{theorem}

{\sc Proof.} According to proposition \ref{PropClosedScalarization} we have
\[
\forall x\in X \colon g\of{x}  = \bigcap\limits_{z^*\in C^-\bs\cb{0}}
    \cb{z\in Z \colon \of{\cl \vp_{g, z^*}^\triup}\of{x} \leq -z^*\of{z}},
\]
and $\vp_{g, z^*}^\triup$ is convex for all $z^*\in C^-\bs\cb{0}$. Applying theorem
\ref{ThmBiconjugationScalar} to $\cl \vp_{g, z^*}^\triup$ we obtain
\begin{align*}
 \forall x\in X \colon g\of{x}
    & = \bigcap\limits_{z^*\in C^-\bs\cb{0}} \cb{z\in Z \colon \of{\vp_{g, z^*}^\triup}^{**}\of{x} \leq -z^*\of{z}}\\
    & = \bigcap\limits_{\substack{(\xi,r)\in X^\triup\times \R\\z^*\in C^-\bs\cb{0}}} \cb{z\in Z \colon \of{\xi_r(x)\idif \of{\vp_{g, z^*}^\triup}^*\of{\xi,r}} \leq -z^*\of{z}}.
\end{align*}
On the other hand, by formula \eqref{EqScalarizationOfConjugate} for all $x\in X$ and
$(\xi,r,z^*)\in X^\triup\times \R\times C^-\setminus\cb{0}$ it holds
\begin{align*}
S^\triup_{(\xi,r,z^*)}(x)\idif g^*(\xi,r)
    = \cb{ z\in Z\colon S^\triup_{(\xi,r,z^*)}(x)\supseteq g^*(\xi,r)+z }\\
    = \cb{ z\in Z\colon \xi_r(x)\leq \of{\vp_{g, z^*}^\triup}^*\of{\xi,r}\isum(-z^*(z)) } .
\end{align*}
Thus, in view of formula \eqref{BiconjugateSetValued3} the statement is proven. \pend

\begin{remark}
Defining the closed convex hull $\cl\co g$ of an arbitrary function $g \colon X \to
Q^\triup$ by $\gr\of{\cl\co g} = \cl\co\of{\gr g}$ we can draw the same conclusion as in
the scalar case since, as in the scalar case, the conjugate of $g$ coincides with the one
of $\cl \co g$ (see e.g. \eqref{EqScalarizationOfConjugate}), namely $g^{**} = \cl \co
g$.
\end{remark}

\section{Appendix}

\subsection{Power sets of linear spaces}

Let $Z$ be a linear space and $\mathcal{P}\of{Z}$ the set of all subsets of $Z$ including
the empty set. The usual Minkowski addition of two sets $A, B \subseteq Z$
\[
A + B = \cb{a + b \colon a \in A, \, b \in B}
\]
is extended to $\mathcal{P}\of{Z}$ by setting $A + B = \emptyset$ if $A = \emptyset$ or
$B = \emptyset$, or both. If $B = \cb{z}$ is a singleton, we abbreviate $A + B = A +
\cb{z} = A + z$.

The point-wise multiplication of a set $A \subseteq Z$ by a non-negative $t \in \R$
\[
t \cdot A = \cb{ta \colon a \in A}
\]
is extended to $\mathcal{P}\of{Z}$ by setting $t \cdot \emptyset = \emptyset$ for all $t
> 0$, and finally $0 \cdot \emptyset = \cb{0}$ where the first $0$ is in $\R$, the second in
$Z$. Finally, we write $A - B$ for $A + \of{-B}$.

\subsection{Conlinear spaces}

The following definition is taken from \cite{Hamel05Habil} where references and more
material about structural properties of conlinear spaces can be found.

\begin{definition}
\label{DefConlinearSpace} A nonempty set $W$ together with two algebraic operations $+
\colon W \times W \to W$ and $\cdot \colon \R_+ \times W \to W$ is called a conlinear
space provided that
\\
(C1) $\of{W, +}$ is a commutative monoid with neutral element $\theta$,
\\
(C2) (i) $\forall w_1, w_2 \in W$, $\forall r \in \R_+$: $r \cdot \of{w_1 + w_2} = r
\cdot w_1 + r \cdot w_2$, (ii) $\forall w \in W$, $\forall r, s \in \R_+$: $s \cdot \of{r
\cdot w} = \of{rs} \cdot w$, (iii) $\forall w \in W$: $1 \cdot w = w$, (iv) $0 \cdot
\theta = \theta$.

An element $w \in W$ is called a convex element of the conlinear space $W$ if
\[
\forall s, t \geq 0 \colon \of{s+t} \cdot w = s \cdot w + t \cdot w.
\]

A conlinear space $\of{W, +, \cdot}$ together with a partial order $\preceq$ on $W$ (a
reflexive, antisymmetric, transitive relation) is called ordered conlinear space provided
that (iv) $w, w_1, w_2 \in W$, $w_1 \preceq w_2$ imply $w_1 + w \preceq w_2 + w$, (v)
$w_1, w_2 \in W$, $w_1 \preceq w_2$, $r \in \R_+$ imply $r \cdot w_1 \preceq r\cdot w_2$.
\end{definition}

\end{document}